\documentclass[preprint]{elsarticle}
\usepackage{amssymb}
\usepackage{amsmath}
\usepackage{graphicx}
\usepackage{tikz}

\usepackage{fancyhdr}
\newcommand\beq{\begin{equation}}
\newcommand\eeq{\end{equation}}

\journal{Journal of Computational and Applied Mathematics}

\begin{document}

\begin{frontmatter}

\title{A multi-level preconditioned Krylov method for the efficient solution of algebraic tomographic reconstruction problems}

\author[address1]{Siegfried Cools\corref{mycorrespondingauthor}}
\cortext[mycorrespondingauthor]{Corresponding author}
\ead{siegfried.cools@uantwerp.be}

\author[address2]{Pieter Ghysels}
\author[address3]{Wim van Aarle}
\author[address3]{Jan Sijbers}
\author[address1]{Wim Vanroose}

\address[address1]{Applied Mathematics Group, University of Antwerp, Middelheimlaan 1, 2020 Antwerp, BE}
\address[address2]{Future Technologies Group, Computational Research Division, Lawrence Berkeley National Laboratory, 1 Cyclotron Road, Berkeley, CA 94720, USA}
\address[address3]{iMinds-Vision Lab, University of Antwerp, Universiteitsplein 1, 2610 Wilrijk, BE \vspace{-0.5cm}}

\begin{abstract}
{\footnotesize
Classical iterative methods for tomographic reconstruction include the class of Algebraic Reconstruction Techniques (ART). Convergence of these stationary linear iterative methods is however notably slow.  
In this paper we propose the use of Krylov solvers for tomographic linear inversion problems. These advanced iterative methods feature fast convergence at the expense of a higher computational cost per iteration, causing them to be generally uncompetitive without the inclusion of a suitable preconditioner. Combining elements from standard multigrid (MG) solvers and the theory of wavelets, a novel wavelet-based multi-level (WMG) preconditioner is introduced, which is shown to significantly speed-up Krylov convergence. 
The performance of the WMG-preconditioned Krylov method is analyzed through a spectral analysis, and the approach is compared to existing methods like the classical Simultaneous Iterative Reconstruction Technique (SIRT) and unpreconditioned Krylov methods on a 2D tomographic benchmark problem. 
Numerical experiments are promising, showing the method to be competitive with the classical Algebraic Reconstruction Techniques in terms of convergence speed and overall performance (CPU time) as well as precision of the reconstruction.
}
\end{abstract}


\end{frontmatter}

\section{Introduction}

Computed Tomography (CT) is a powerful imaging technique that allows non-destructive visualization of the interior of physical objects.  Besides its common use in medical applications \cite{arridge1999optical}, tomography is also widely applicable in fields such as biomedical research, materials science, metrology, etc.
In all applications, a certain imaging source (e.g.~an X-ray source) and an imaging detector (e.g.~X-ray detector) are used to acquire two-dimensional projection images of the object from different directions. A three-dimensional virtual reconstruction can then be computed using one of the many reconstruction techniques that can be found in the literature. In practice, the most commonly used analytical methods for CT are Filtered Backprojection (FBP) and its cone-beam variant Feldkamp-Davis-Kress (FDK).  These methods make use of various analytical properties of the projection geometries to compute the reconstructed object at a low computational cost. A major drawback of analytical methods is their inflexibility to different experimental setups and their inability to include reconstruction constraints which can be used to exploit possible prior information about the object. 

Iterative Algebraic Reconstruction Techniques (ART) form an interesting alternative to the aforementioned analytical methods.  Here, the reconstruction problem is described as the solving of a system of linear equations. The Simultaneous Iterative Reconstruction Technique (SIRT) is a straightforward method that has been extensively studied in the literature, see \cite{Gregor2008} and the references therein.  Another general class of algebraic solution methods are the Krylov solvers such as CGLS, GMRES, BiCGStab, etc., an overview of which can be found in \cite{simoncini2007recent}.  Alternatively, one can resort to more powerful techniques that apply additional constraints to the reconstruction, which can lead to improved accuracy, especially when fewer projection images are available (i.e.~scans with a lower radiation dose). Total variation minimization approaches such as FISTA \cite{Beck2009}, for example, assume that the variation between neighbouring pixels is low inside a homogeneous object.  Discrete tomography approaches such as DART \cite{Batenburg2011} improve the reconstruction quality by limiting the number of grey level values that can be present in the reconstructed image.

While iterative methods for tomography have become widely accepted in the scientific community, practical applications have not yet adopted these techniques \cite{Pan2009}, mostly due to the variable computational cost and storage requirements of the iterative process (contrary to the fixed costs of analytical methods based on FFT-type algorithms).  The development of efficient new iterative solvers is therefore crucial.  This efficiency can be accomplished in two ways.  Firstly, the computation time of each iteration can be reduced by optimally exploiting parallelism of the projection and backprojection operators with the use of modern hardware accelerated computer architectures such as NVIDIA GPU's \cite{Palenstijn2011} or the Intel Xeon Phi \cite{vanAarle2012memory}. Secondly, a solver with a fast convergence rate, requiring only a limited number of iterations should be used. Additionally, the convergence rate of the ideal solver should not depend on the problem size. 

In this work, an approach that fits into the second category will be introduced for non-constrained iterative reconstruction.  By analyzing the spectral properties of the standard SIRT method, it will be shown that the convergence of classical algebraic reconstruction techniques (stationary iterative schemes) is notably slow. As it appears, the alleged smoothing property does not hold in the case of tomographic reconstruction problems. Krylov methods prove to be more efficient, yet are generally more expensive in terms of memory and computation cost.  Therefore, when using Krylov methods, it is mandatory to define an efficient preconditioner, which allows faster convergence.  This approach is very common in a wide range of PDE-type problems, yet is still fairly new for tomographic reconstruction. Related work in the setting of tomographic reconstruction includes the research on multilevel image reconstruction by McCormick et al.~\cite{mccormick1993multigrid, henson1996multilevel}, and more recently the work done on multigrid methods for tomographic reconstruction by Webb et al.~\cite{oh2005general} and R\"ude et al.~\cite{kostler2006towards}.

Originally introduced as a theoretical tool by Fedorenko in 1964 \cite{fedorenko1964speed} and later adopted as a solution method by Brandt in 1977 \cite{brandt1977multi}, multigrid (MG) solvers are commonly used as efficient and low-cost Krylov preconditioners for high-dimensional problems in the PDE literature, see e.g.~\cite{elman2002multigrid,erlangga2009algebraic}. One of the key concepts of the multigrid scheme is the representation of the original fine grid reconstruction problem on a coarser scale resolution, where the problem is computationally cheaper to solve. However, we show that the standard multigrid approach \cite{briggs2000multigrid,stüben1982multigrid,trottenberg2001multigrid,hackbusch1985multi} does not act as an efficient preconditioner for algebraic tomographic systems. Indeed, the ineffectiveness of the smoother in eliminating the oscillatory modes causes the key complementary action of smoother and coarse grid correction to fail, resulting in an inefficient multigrid scheme for algebraic tomographic reconstruction problems.

In this work a new wavelet-based multigrid (WMG) preconditioner is introduced, which is more suited for tomographic reconstruction. The proposed method combines elements from standard multigrid with the theory of wavelets, and shows some similarities to the work on wavelet-based multiresolution tomographic reconstruction in \cite{delaney1995multiresolution} and \cite{bhatia1996wavelet}. Additionally, the main advantage of the proposed method, i.e.\,projection of the large fine-scale system onto smaller, easy-to-solve subproblems, resembles key features of the Hierarchical Basis Multigrid Method (HBMM) \cite{bank1988hierarchical,bank1996hierarchical}. It is shown through an eigenvalue analysis that WMG-preconditioning significantly increases Krylov convergence speed, which is confirmed by various numerical experiments. Additionally, we show that the WMG-preconditioned Krylov solver allows for an accuracy which is generally unobtainable by classical SIRT reconstruction. The numerical results presented in this work show promise, validating the proposed WMG scheme as an efficient Krylov preconditioning technique for algebraic tomographic reconstruction.

The paper is structured as follows. In Section 2 the classical SIRT and MG-Krylov solvers for iterative tomographic reconstruction are reviewed and analyzed. Section 3 introduces a novel preconditioning approach to account for the defects of the MG preconditioner, which greatly improves convergence speed of the BiCGStab Krylov solver. In Section 4, a series of experimental simulations is presented to validate our contribution. Ultimately, Section 5 concludes this work with an overview of the main results in this paper and a discussion on possible future research options.

\section{Notation and key concepts of tomographic reconstruction}

\subsection{Algebraic tomographic reconstruction}

Consider a data vector $b \in \mathbb{R}^{M}$, with $M = m \times n$, where $m$ is the number of projection angles and $n$ is the number of beams. We assume that the number of pixels in every spatial dimension equals $n$, such that the data is reconstructed on a 2D $n \times n$ grid. We denote the total number of pixels in the image by $N = n \times n$. Algebraic reconstruction methods consider tomographic reconstruction as the problem of solving the linear system of equations
\beq \label{eq:sys}
W x = b,
\eeq
where $x = (x_j) \in \mathbb{R}^{N}$ are the unknown attenuation values on the
grid in the image domain, which represent the object of interest, and $b = (b_i) \in \mathbb{R}^{M}$ are the measured projection
values for each beam and under each angle. The matrix $W = (w_{ij}) \in \mathbb{R}^{M \times N}$ is a linear projection operator that maps an image in the object or reconstruction domain onto the projection domain, see Figure \ref{fig:pixelsfinecoa}. It can be computed or approximated in a variety of ways. In this work, Joseph's projection kernel \cite{Joseph1982} will be used. Note that the matrix $W$ is generally large and sparse, featuring $\mathcal{O}(\sqrt{N})$ non-zero elements per column.

For solving purposes, system (\ref{eq:sys}) is frequently rewritten as the equivalent system of normal equations (NE)
\beq \label{eq:sys_ne}
W^T W x = W^T b,
\eeq
replacing the (possibly underdetermined) rectangular system (\ref{eq:sys}) by a system with a square symmetric matrix $W^T W \in \mathbb{R}^{N\times N}$. The normal form is commonly used when solving the above system for $x$ using Algebraic Reconstruction Techniques or Krylov methods as described in Section 2.2. In view of efficient implementation, the matrix multiplication $W^T W$ is never computed explicitly, as this would result in a dense matrix of $\mathcal(O)(N)$ non-zero elements per column. Instead, the application of $W^T W$ to a vector $x$ is computed as two sparse matrix-vector products (SpMV): $W^T (W x)$.

For convenience of the analysis we consider an ideal experiment without the incorporation of noise throughout the following sections. Additionally, we assume that a sufficient number of projections is given, such that (\ref{eq:sys}) has a unique solution, i.e.~we effectively assume that problem (\ref{eq:sys}) is well-posed. We refer to Section 4.2 for the more realistic case of noisy and/or low-data problems, where regularization is introduced to account for the ill-posedness of the problem.

We first discuss some well-known classes of iterative methods for the tomographic system (\ref{eq:sys})-(\ref{eq:sys_ne}). It will be shown that standard Algebraic Reconstruction Techniques (ART) fail to efficiently solve the problem. Krylov methods prove to be more efficient, yet require a suitable preconditioner to allow for fast convergence. We show that classical multigrid preconditioning is not efficient in the context of tomography, which motivates the introduction of a new multi-level type preconditioner in Section 3.

\subsection{Classical algebraic reconstruction techniques} \label{sec:classical_art}

Algebraic reconstruction algorithms are among the current-day state-of-the-art methods for solving tomographic systems. The results presented in this section are well-known in the literature, see standard works on the basic principles of computerized tomography \cite{natterer1986computerized,kak2001principles} and the references therein. However, the presented insights provide a strong motivation for the construction of the WMG preconditioner in Section 3.

\subsubsection{SIRT method.} The Simultaneous Iterative Reconstruction Technique (SIRT) is a basic stationary iterative scheme for the solution of linear systems of equations, which aims at solving (\ref{eq:sys}) iteratively using a basic residual estimation scheme. Consider the scaled system of normal equations, equivalent to (\ref{eq:sys}),
\beq \label{eq:sys_sirt}
C W^T R W x = C W^T R b,
\eeq
where $R = (r_{ij}) \in \mathbb{R}^{M \times M}$ is a diagonal matrix of the inverse row sums of $W$,
\beq
r_{ii} = \left(\sum_{j=1}^{N} w_{ij}\right)^{-1} \textrm{for} \quad i = 1,\ldots,M; \qquad r_{ij} = 0 \quad \textrm{for} \quad i \neq j.
\eeq
Likewise, $C = (c_{ij}) \in \mathbb{R}^{N\times N}$ is a diagonal matrix of the inverse column sums of $W$,
\beq
c_{jj} = \left(\sum_{i=1}^{M} w_{ij}\right)^{-1} \textrm{for} \quad j = 1,\ldots,N; \qquad c_{ij} = 0 \quad \textrm{for} \quad i \neq j.
\eeq
The scaled system \eqref{eq:sys_sirt} was proposed by Gregor and Benson in \cite{Gregor2008}, where it was shown that this scaling is mandatory to ensure the stability of the SIRT scheme. The SIRT iteration scheme can be written recursively as a stationary iteration
\beq \label{eq:sirt_it}
x^{(k+1)} = x^{(k)} + r^{(k)} =  x^{(k)} + C W^T R (b - W x^{(k)}), \qquad k \in \{1,2,\ldots\}.
\eeq
Note that the SIRT scheme allows for a matrix-free implementation; in practice the operator $C W^T R W$ is never formed explicitly, but its application is typically implemented as a series of SpMVs or a matrix-free projection simulator. Key features of the SIRT method are its low storage cost (only the current guess $x^{(k)}$ and the residual $r^{(k)} = b - W x^{(k)}$ need to be stored) and its relatively low computational cost per iteration (only two SpMV operations, which can be easily parallelized). However, as will be discussed in Section 2.2.2, the convergence rate of the SIRT scheme for tomographic reconstruction is very slow. A large number of iterations is typically required, implying long overall computational times. Nevertheless, due to its straightforward implementation SIRT is commonly used throughout the scientific literature and in practical tomographic implementations. In this work, the SIRT solver will act as the benchmark algebraic solution scheme for tomographic reconstruction.

\begin{figure}[t]
\begin{center}
\includegraphics[width=0.29\textwidth]{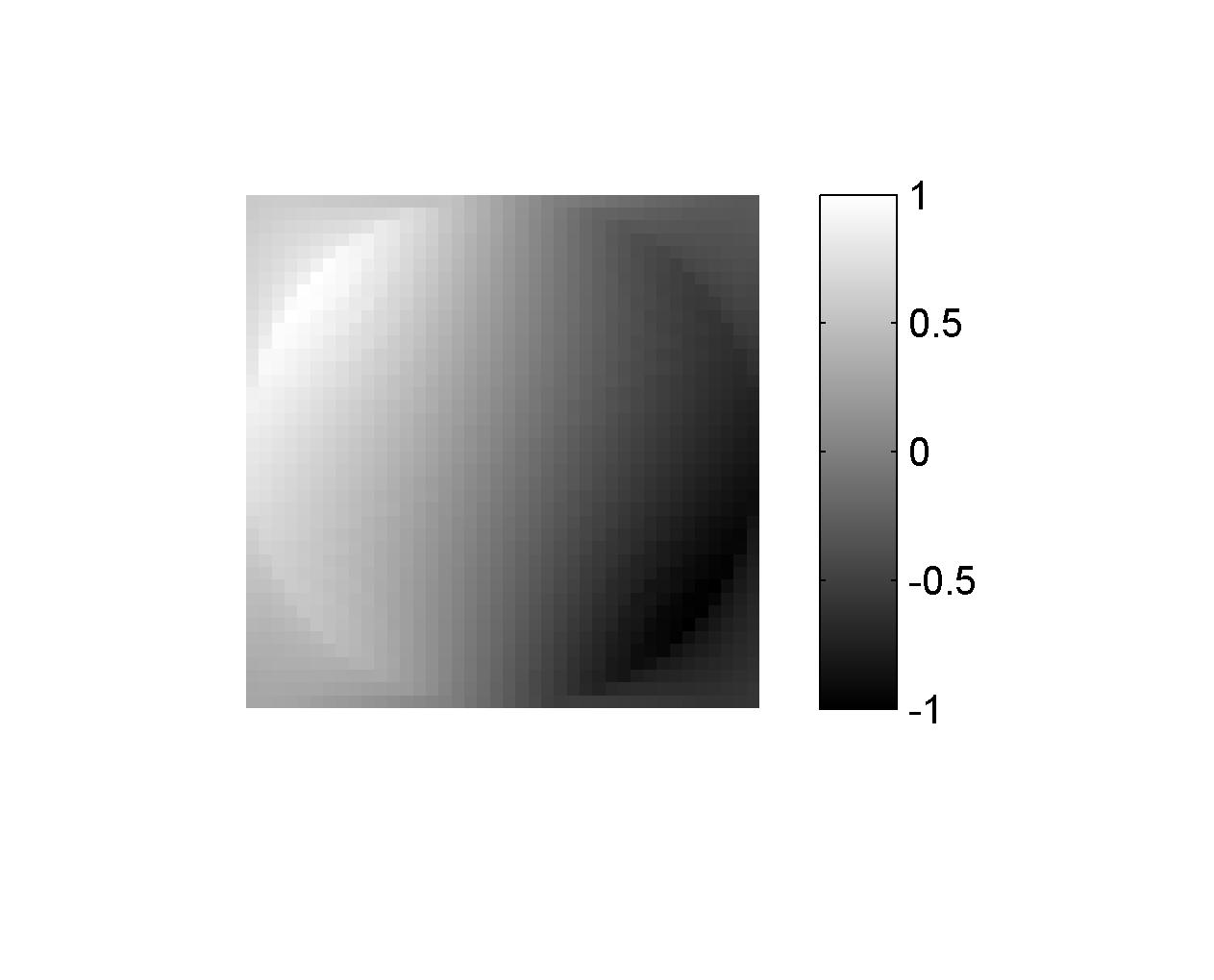} \hspace{-0.9cm}
\includegraphics[width=0.29\textwidth]{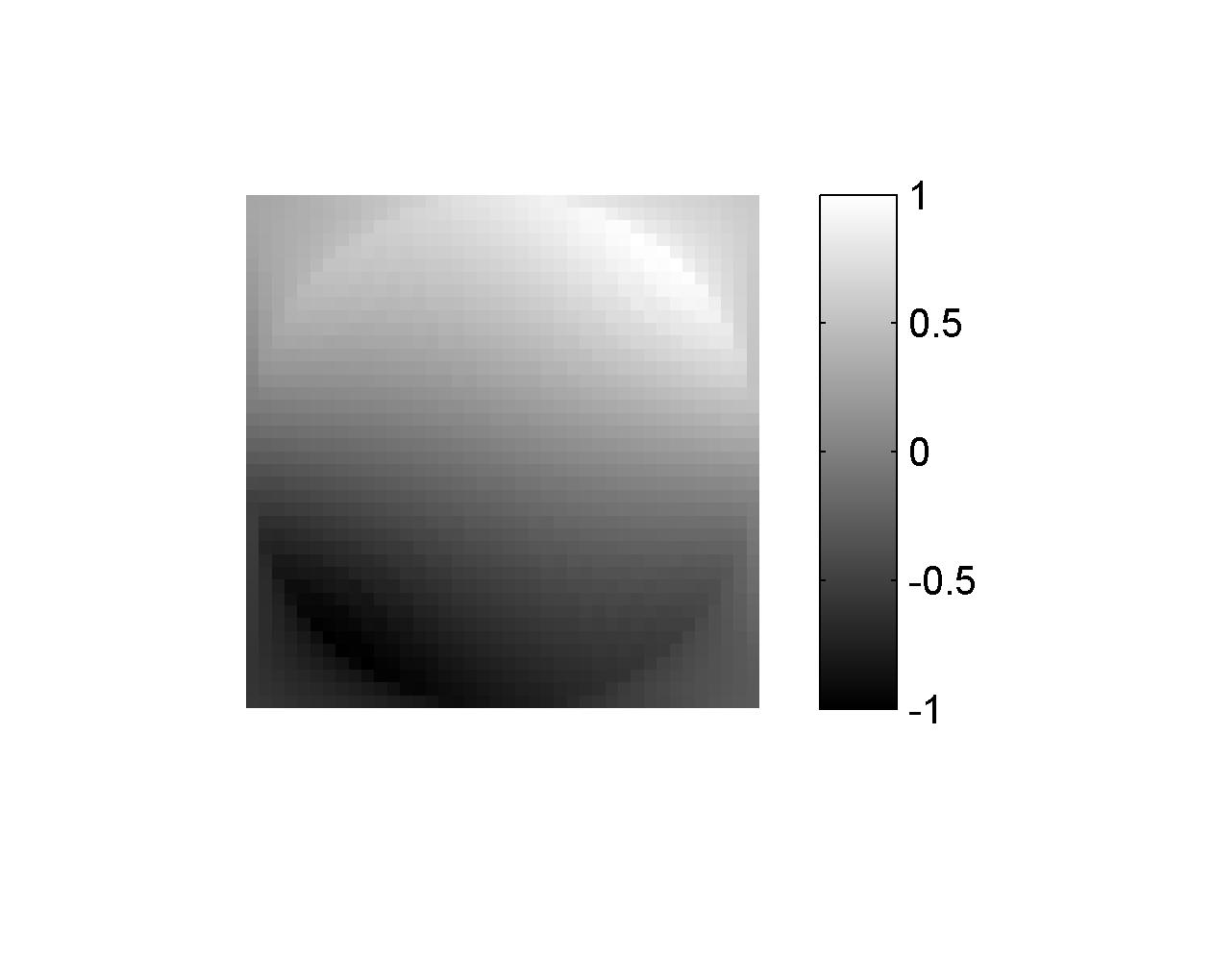} \hspace{-0.9cm}
\includegraphics[width=0.29\textwidth]{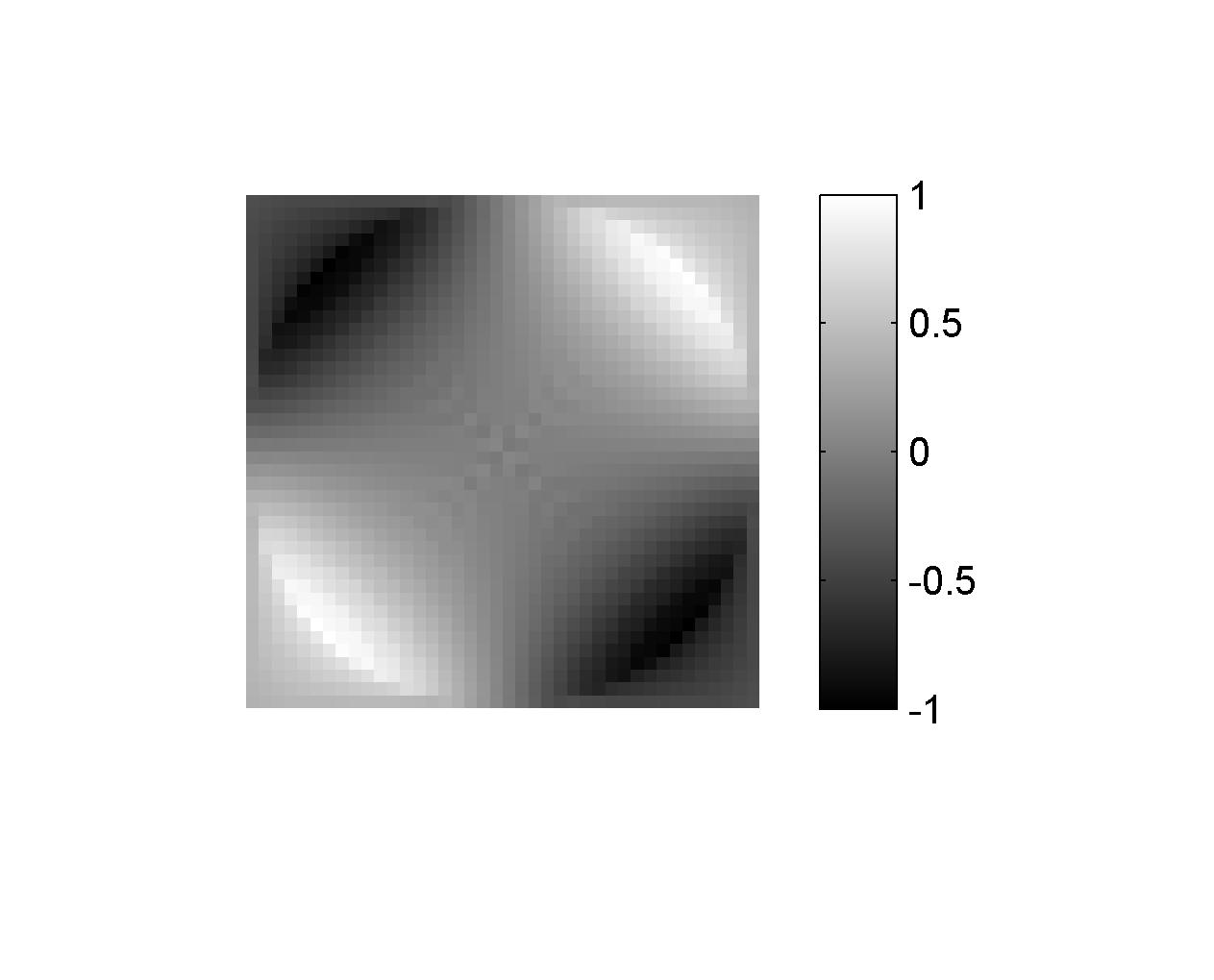} \hspace{-0.9cm}
\includegraphics[width=0.29\textwidth]{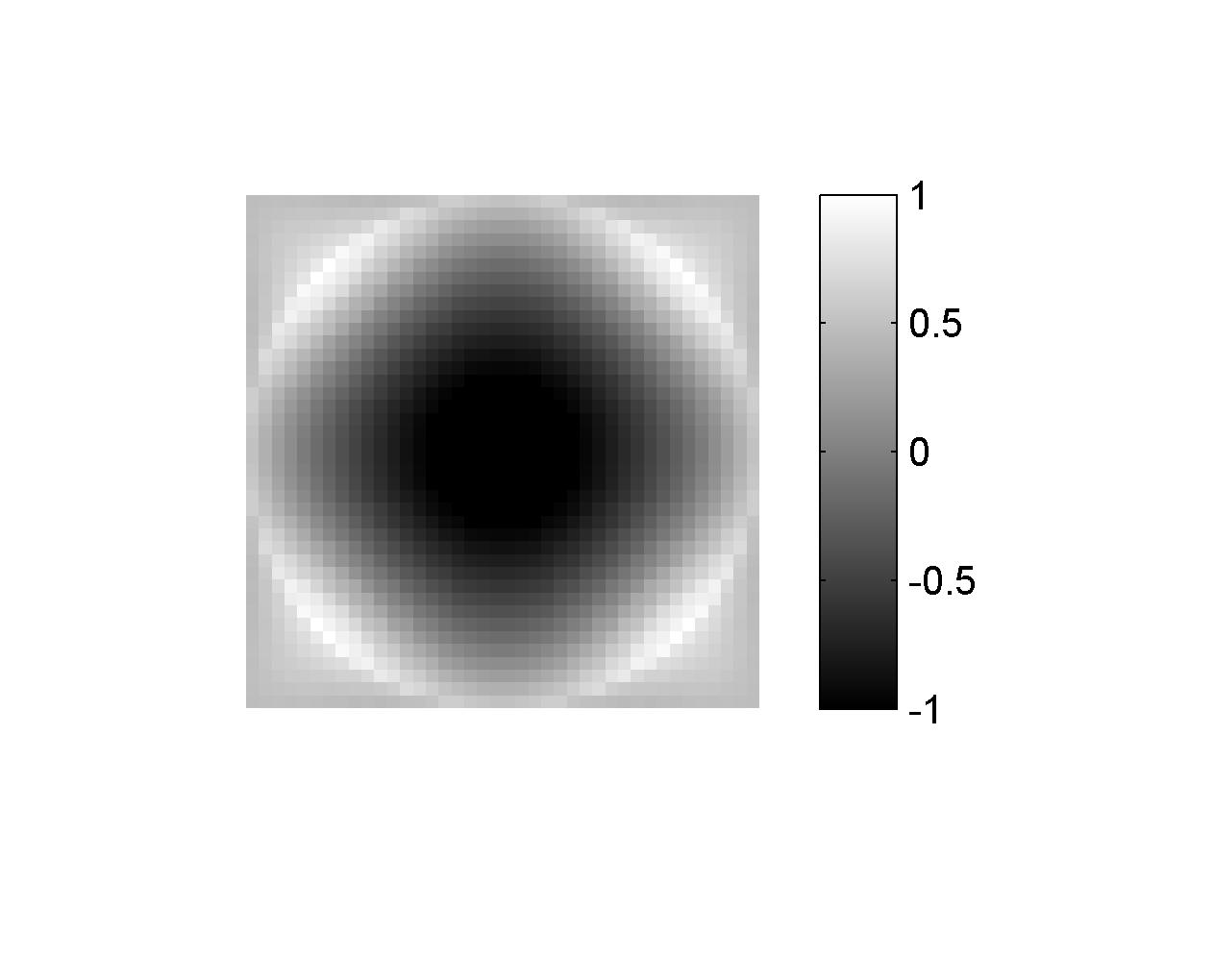} \\ \vspace{-0.7cm}
\includegraphics[width=0.29\textwidth]{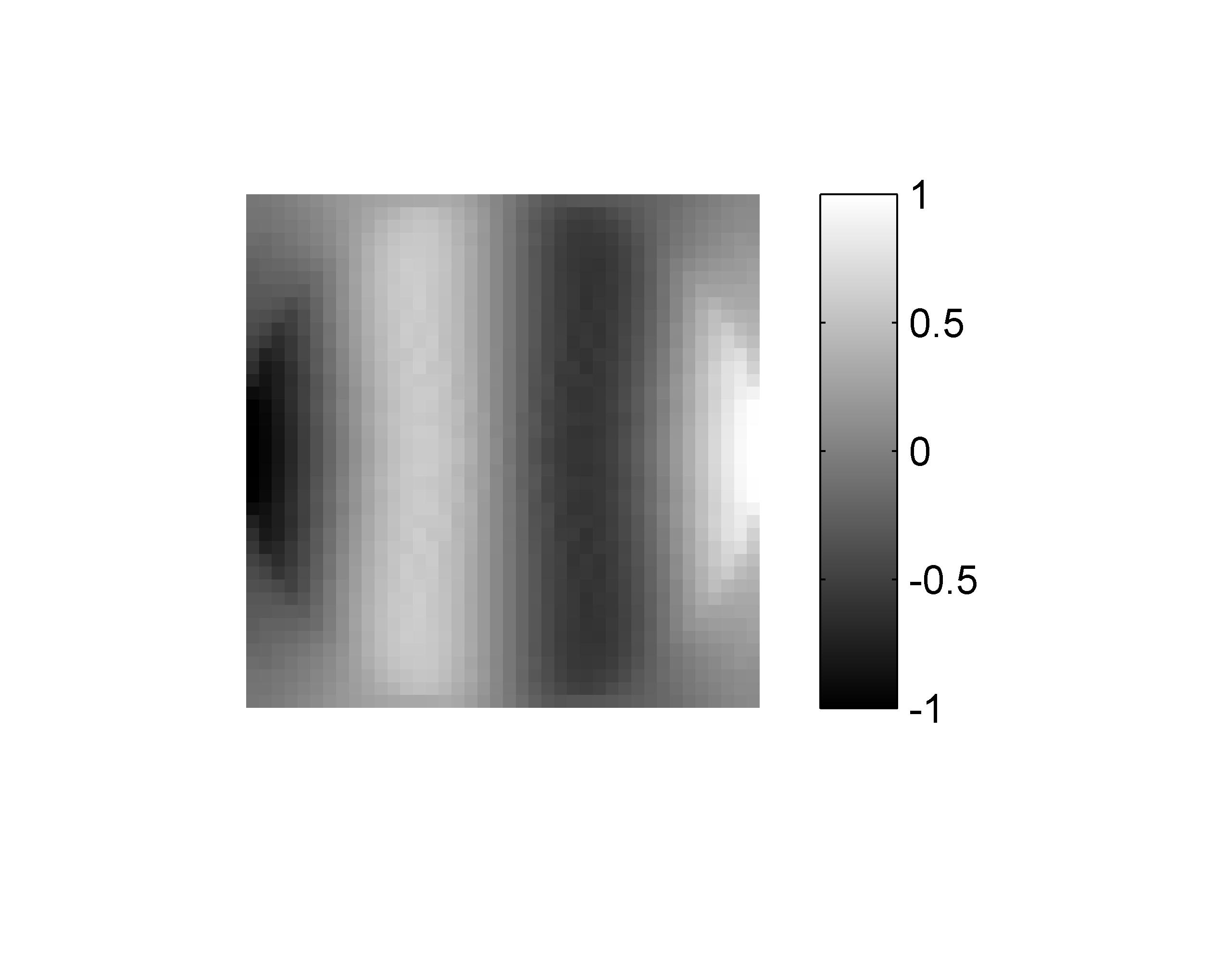} \hspace{-0.9cm}
\includegraphics[width=0.29\textwidth]{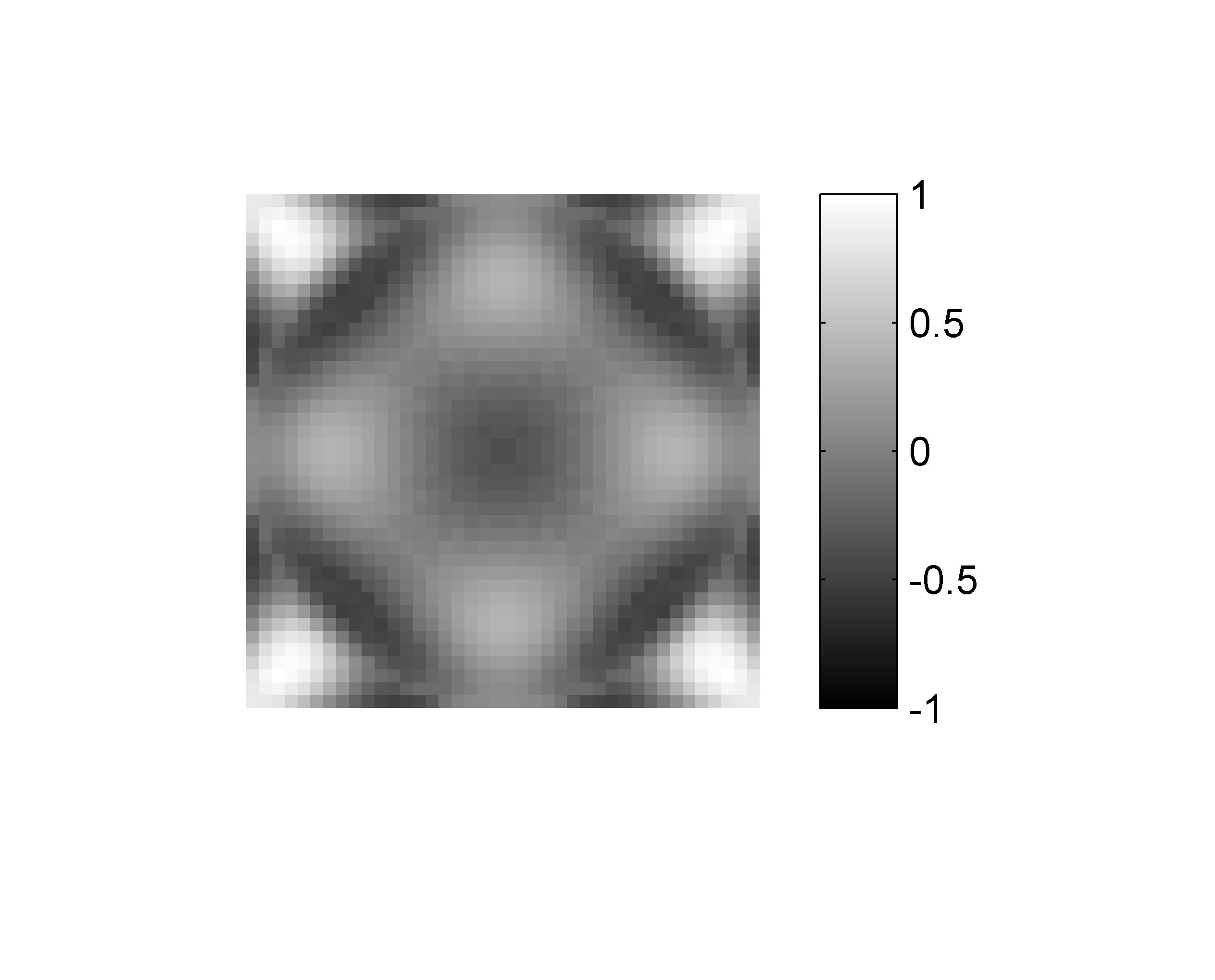} \hspace{-0.9cm}
\includegraphics[width=0.29\textwidth]{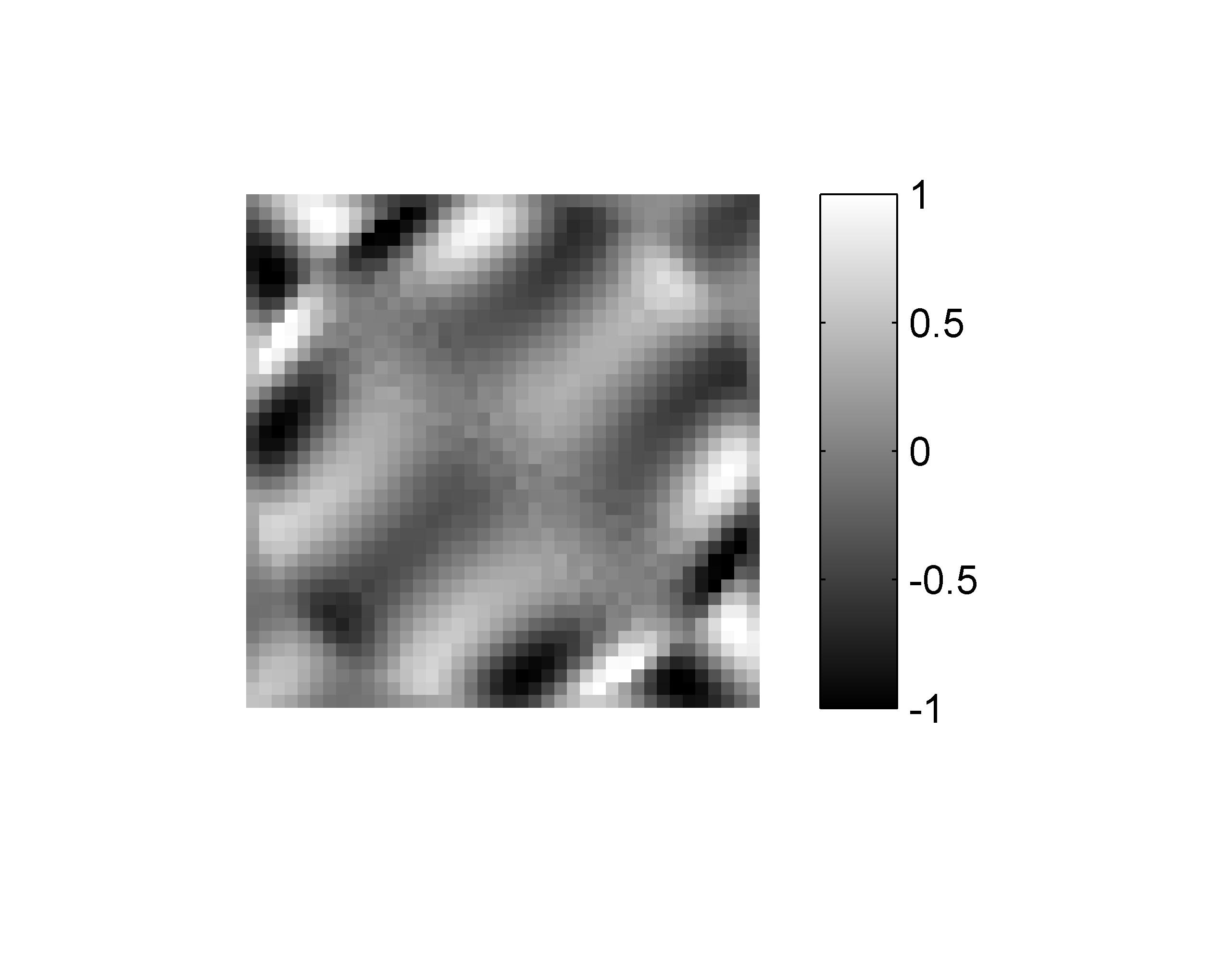} \hspace{-0.9cm}
\includegraphics[width=0.29\textwidth]{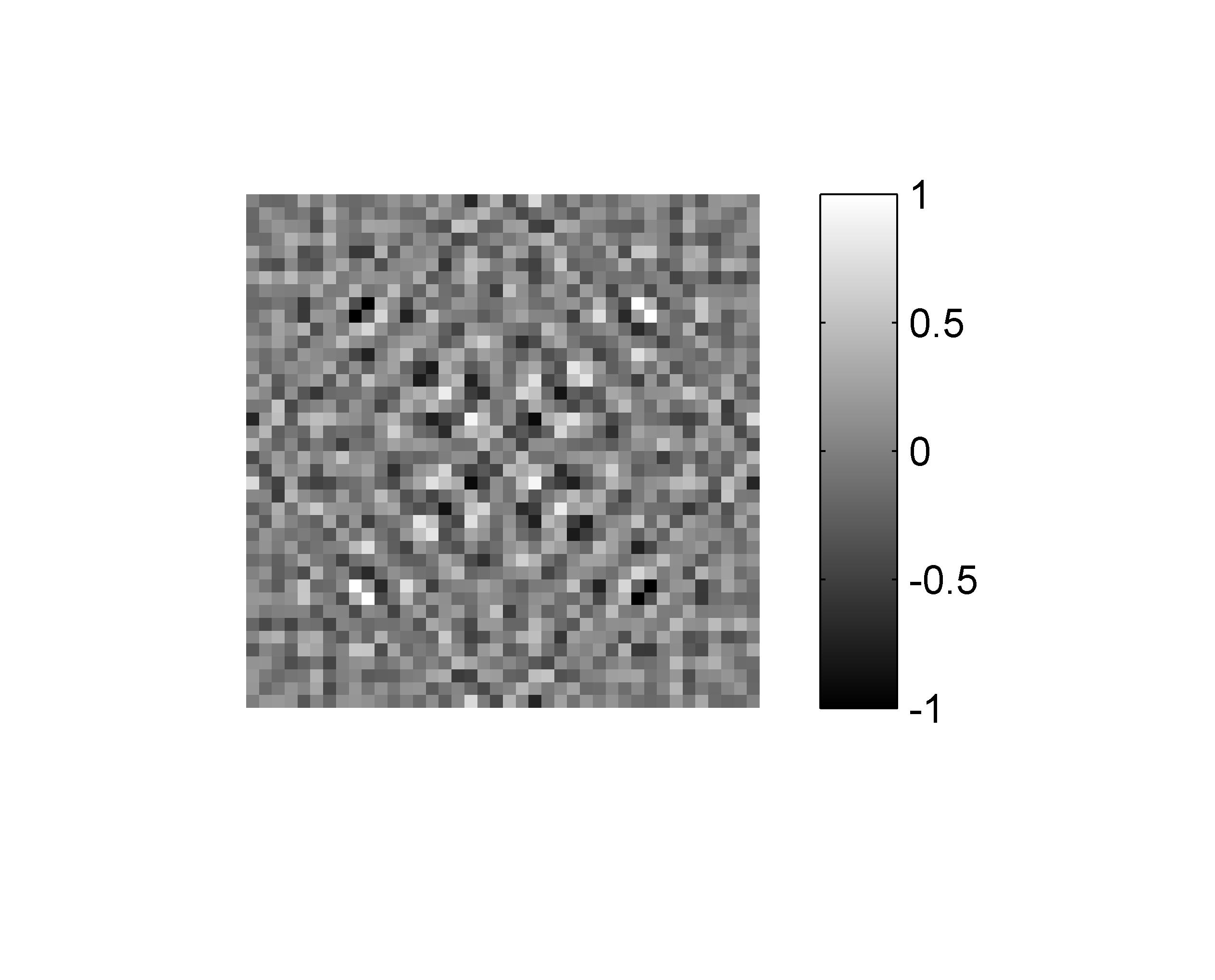} \\ \vspace{-0.5cm}
\caption{A selection of eigenmodes $\varphi_i$ of the SIRT iteration matrix $S$ (\ref{eq:ssirt}) with $N = 40 \times 40$ and $M = 100 \times 40$. First row (f.l.t.r.): very smooth modes $\varphi_2$, $\varphi_3$, $\varphi_4$ and $\varphi_5$. Second row (f.l.t.r.): smooth modes $\varphi_{10}, \varphi_{20}$ and $\varphi_{50}$, and the oscillatory mode $\varphi_{800}$. Corresponding eigenvalues: $\lambda^S_2 = 0.52$, $\lambda^S_3 = 0.52$, $\lambda^S_4 = 0.66$, $\lambda^S_4 = 0.70$, $\lambda^S_{10} = 0.80$, $\lambda^S_{20} = 0.86$, $\lambda^S_{50} = 0.92$, and $\lambda^S_{800} = 0.98$.
}
\label{fig:eig_vec_sirt}
\end{center}
\end{figure}

\subsubsection{Spectral analysis of SIRT.} The origin of the slow SIRT convergence can be found by analyzing the spectrum of the iteration matrix. From the recursion (\ref{eq:sirt_it}), the SIRT method can be interpreted as a basic stationary iterative scheme (cf.~Richardson iteration, weighted Jacobi or Gauss-Seidel), with iteration matrix $S$ given by
\beq \label{eq:ssirt}
S = I - C W^T R W.
\eeq
Consequently, the error in every iteration $k \in \{1,2,\ldots\}$ can be written as
\beq
e^{(k)} = S e^{(k-1)} = S^{k} e^{(0)}, \qquad k \in \{1,2,\ldots\}.
\eeq
Writing the error $e^{(k)}$ as a linear combination of the eigenmodes $\varphi_i$  $(i=1,\ldots,N)$ of $C W^T R W$, which notably are also eigenmodes of $S$, we have
\beq
e^{(k)} = \sum_{i=1}^{N} \alpha^{(k)}_i \varphi_i, \qquad k \in \{1,2,\ldots\}.
\eeq
The eigenmodes of $S$ are interpreted as basis functions for the error after each iteration, and the eigenvalues of $S$, $\lambda^S_i$  $(i=1,\ldots,N)$, represent the propagation factors for the basis functions in the error, as 
\beq \label{eq:sirt_prop}
e^{(k)} = S e^{(k-1)} = \sum_{i=1}^{N} \lambda^S_i \alpha^{(k)}_i \varphi_i, \qquad k \in \{1,2,\ldots\}.
\eeq
This implies that the reduction of the error components $\varphi_i$ in each SIRT iteration is governed by the eigenvalues of the iteration matrix $S$. 
Figure \ref{fig:eig_vec_sirt} shows a number of eigenmodes $\varphi_i$ of the SIRT iteration matrix $S$ for a volume size of $N=40\times40$ with 100 equiangular parallel beam projections over $180^\circ$ of 40 rays each. The eigenmodes are ordered according to the magnitude of their corresponding eigenvalues $\lambda^S_i$, from small (low index) to large (high index). Low-indexed eigenmodes are slow-varying across the numerical domain in both directions, while high-indexed eigenvalues correspond to highly oscillatory modes. For 2D problems, the first quarter of eigenmodes in the spectrum ($\varphi_i$ with $i \leq N/4$) are commonly referred to as smooth modes, while the remaining eigenmodes ($i > N/4$) are oscillatory in one or both spatial directions \cite{briggs2000multigrid}. 

\begin{figure}[t]
\begin{center}
\includegraphics[width=0.48\textwidth]{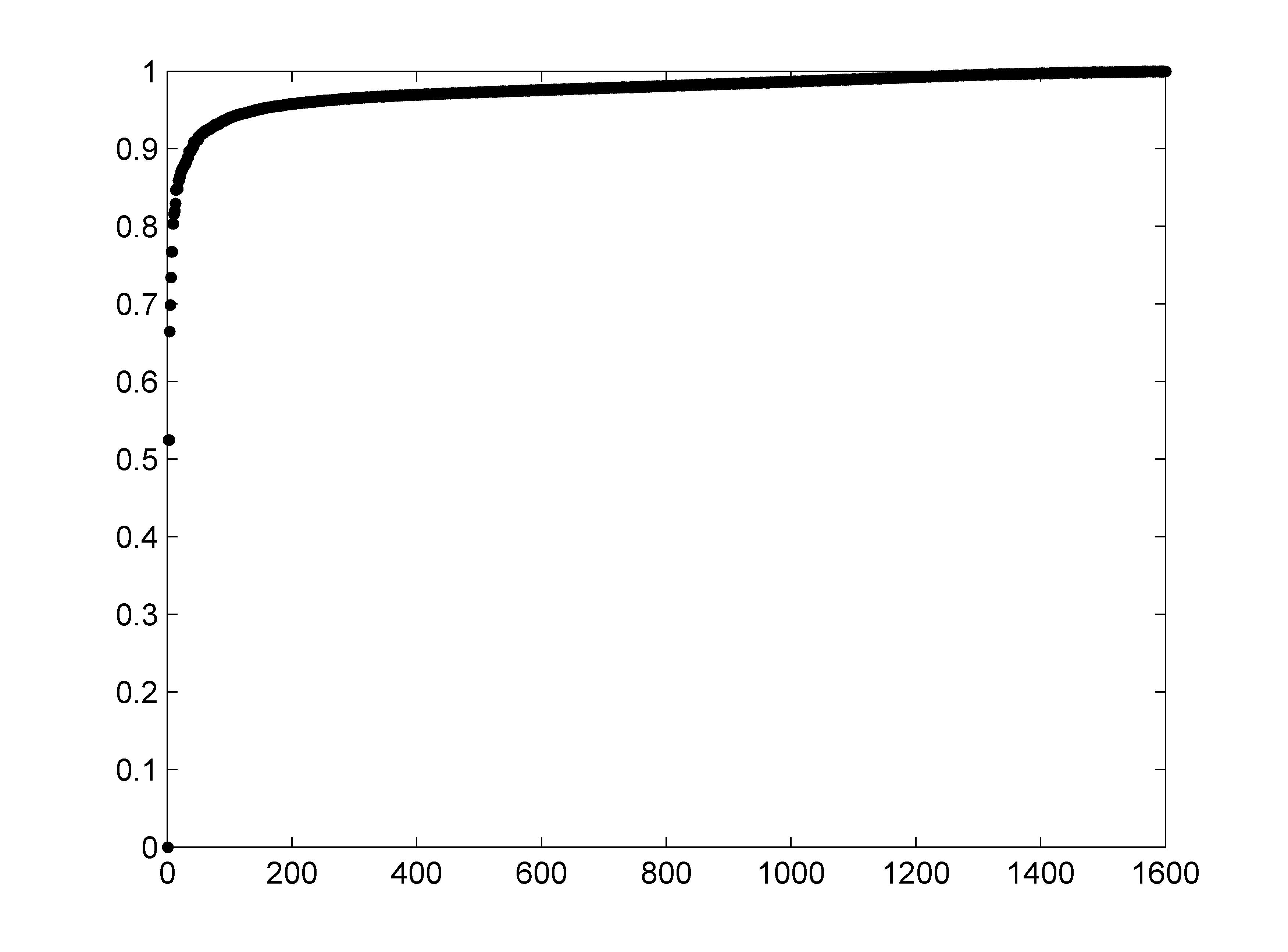} 
\caption{Spectrum $\lambda^S_i$ $(i=1,\ldots,N)$ of the SIRT iteration matrix $S$, see (\ref{eq:ssirt}), with $N = 40 \times 40$ and $M = 100 \times 40$. Eigenvalues act as propagation factors for the error basis functions (eigenmodes), see (\ref{eq:sirt_prop}).}
\label{fig:eig_val_sirt}
\end{center}
\end{figure}

A characteristic property of all stationary linear iterative schemes (alternatively called \emph{relaxation schemes} or \emph{smoothers}) in the setting of Laplace-type PDEs is the damping of oscillatory eigenmodes, while preserving the smooth components in the error \cite{elman2002multigrid,briggs2000multigrid}. This is commonly referred to as the \emph{smoothing property}. However, it appears that this property generally does not hold in the setting of tomographic reconstruction. The numerical eigenvalues $\lambda^S_1, \ldots , \lambda^S_N$ of the SIRT iteration matrix $S$ are shown in Figure \ref{fig:eig_val_sirt}. 
A very limited number of eigenvalues are significantly smaller than one, implying only a small fraction of the error components is damped through successive SIRT iterations.
Furthermore, the smallest eigenvalues correspond exactly to the very \emph{smoothest} eigenmodes of $S$, cf.~Figure \ref{fig:eig_vec_sirt}, which directly contradicts the smoothing property. 

Hence, only a small subset of very smooth eigenmodes are effectively eliminated in every SIRT iteration. Convergence of the SIRT solver -- and thus, by extension, all basic stationary iterative methods -- is therefore notably slow for tomographic reconstruction problems. In addition, these basic stationary iterative methods primarily eliminate the smooth error components, contrary to possessing the smoothing property. These observations provide an important motivation for the use of Krylov methods and the construction of a novel multi-level preconditioner.

\subsection{Krylov methods and classical multigrid preconditioning} 

\subsubsection{Krylov methods.} Primarily used in the solution of high-dimensional PDE's, Krylov methods are less well-known as a class of iterative solvers in the context of tomographic reconstruction. Consider a general linear system of the form 
\beq \label{eq:kryl_sys}
A x = f,
\eeq
where $A \in \mathbb{R}^{N \times N}$ and $f \in \mathbb{R}^N$. Note that the tomographic system (\ref{eq:sys_ne}) is of this form with $A = W^T W$ and $f = W^T b$. In every Krylov iteration, the residual (or some other vector quantifying deviation from the solution) is minimized over the $k$-th Krylov subspace
\beq \label{eq:kryl_space}
\mathcal{K}_k(A,f) = \operatorname{span} \, \{ f, A f, A^2 f, \ldots, A^{k-1} f \}.
\eeq
Many varieties of Krylov solvers may be used to solve system (\ref{eq:kryl_sys}): GMRES, MINRES, CG, BiCG, CGLS, LSQR, etc. In this work we consider BiCGStab \cite{van1992bi} as the primary Krylov solver for system (\ref{eq:sys_ne}), where we assume that $A = W^T W$ and $f = W^T b$ in the above definitions (\ref{eq:kryl_sys})-(\ref{eq:kryl_space}). 

\begin{figure}[t]
\begin{center}
\includegraphics[width=7.5cm]{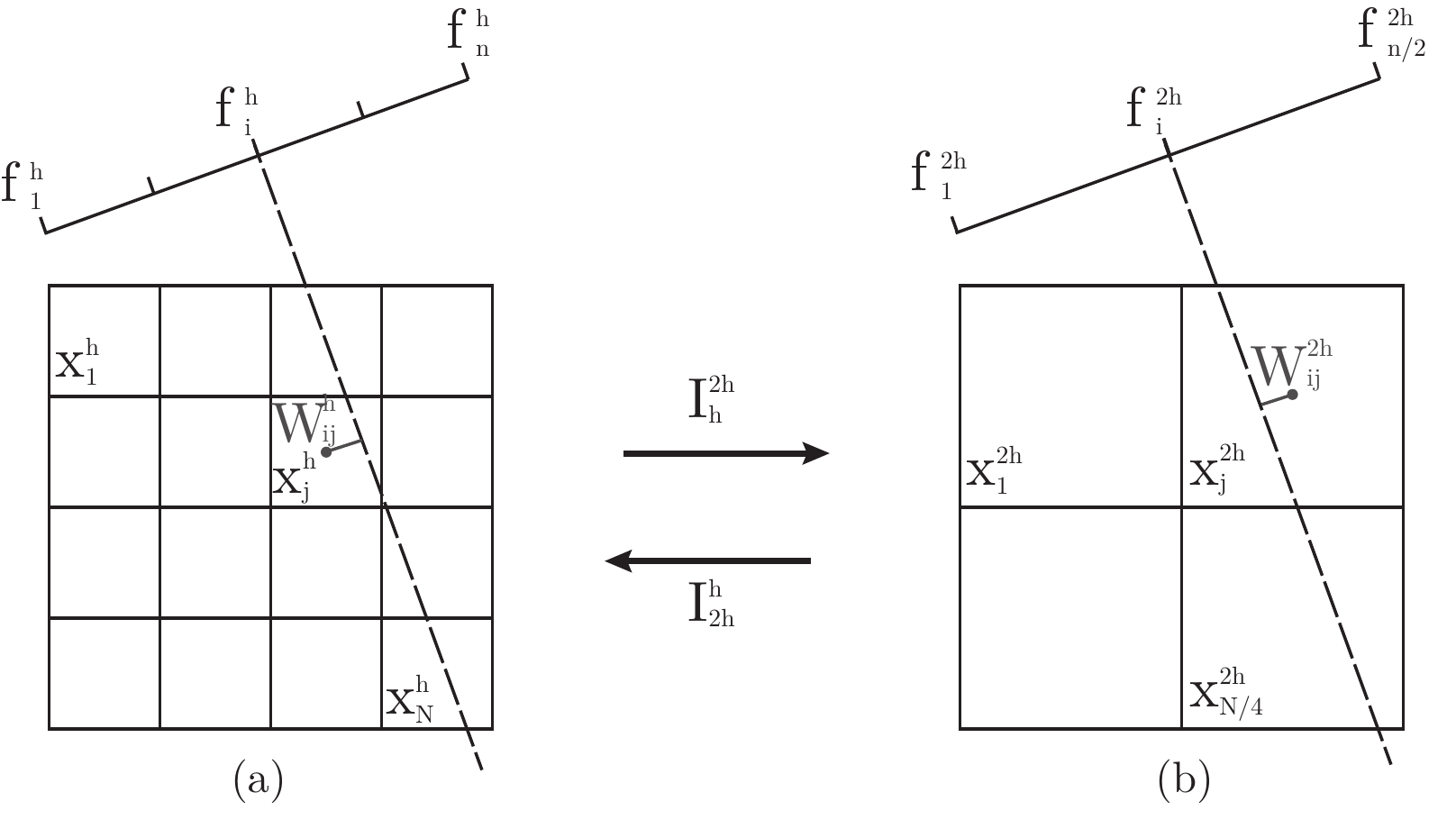}
\caption{Schematic representation of the fine $h$-grid $\Omega^{h}$ (a) and coarse $2h$-grid $\Omega^{2h}$ (b) for 2D multi-level tomographic reconstruction. The projection matrix $W^{2h}$ is redefined to the coarse grid setting.}
\label{fig:pixelsfinecoa}
\end{center}
\end{figure}

While Krylov methods converge considerably faster than most stationary linear iterative schemes, their per-iteration computational cost and storage requirements are generally higher. To keep the number of iterations as low as possible, Krylov methods often incorporate preconditioning by multiplying the system matrix $A$ (from the left or the right) by $M^{-1}$.
For the tomographic reconstruction problem (\ref{eq:sys_ne}), the Krylov method is either applied to the left preconditioned system
\beq \label{eq:sys_prec}
M^{-1} W^T W x = M^{-1} W^T b, 
\eeq
or the right preconditioned system
\beq
W^T W M^{-1} y = W^T b, \quad x = M^{-1} y,
\eeq
for some preconditioning operator $M^{-1} \in \mathbb{R}^{N \times N}$. The application of the operator $M^{-1}$ to a vector $v \in \mathbb{R}^{N}$ is often computed by some numeric scheme (possibly different from the scheme used to compute the SpMV's with $W^T W$), thus avoiding explicit computation of the inverse. Note that the experiments in this paper are based on right preconditioning.

\subsubsection{Multigrid preconditioning.} In this work we propose preconditioning of the system of normal equations (\ref{eq:sys_ne}) by a multi-level type scheme \cite{brandt1977multi}. In recent years, standard multigrid methods \cite{briggs2000multigrid,stüben1982multigrid,trottenberg2001multigrid,hackbusch1985multi} have been broadly used as Krylov preconditioners in various application areas, e.g.~seismic imaging, see \cite{calandra2012two,reps2011analyzing}. 

Let $\Omega^h$ denote a grid with pixel size $h$, and let $A^h$, $x^h$ and $f^h$ respectively denote the system matrix, solution vector and data vector represented on $\Omega^h$. The main idea of multi-level schemes is to represent the original fine grid equation
\beq
A^h x^h = f^h, \qquad x^f \in \Omega^h,
\eeq
on a coarser grid $\Omega^{2h}$, which consists of bigger pixels formed by 2-by-2 blocks of pixels of the original fine resolution grid $\Omega^h$, see Figure \ref{fig:pixelsfinecoa}. To convert data from the fine to coarse grid and vice versa, two intergrid operators are defined: the restriction operator $I_h^{2h}$ and the interpolation operator $I_{2h}^h$. The main advantage of a multi-level approach is that the system is much cheaper to solve numerically on the coarse $\Omega^{2h}$ grid. The coarse grid matrix $A^{2h}$ is either formed by recalculating the attenuation values for the $\Omega^{2h}$ grid explicitly, or through Galerkin coarsening,
\beq \label{eq:galerkin}
A^{2h} = I_{h}^{2h} A^h I_{2h}^h.
\eeq
The restriction and interpolation operators used for 2D tomographic reconstruction in this work are based upon the following one-dimensional intergrid operators:
\beq \label{eq:intops}
I_{h,1D}^{2h} = \frac{1}{\sqrt{2}}
\begin{bmatrix}
1 & 1 & & & & & &\\
& & 1 & 1 & & & & \\
& & & & \cdot &\cdot & &\\
& & & & & & 1 & 1 \\
\end{bmatrix} \in \mathbb{R}^{\frac{n}{2} \times n},
~~
I_{2h,1D}^{h} = (I_{h,1D}^{2h})^T \in \mathbb{R}^{n \times \frac{n}{2}}.
\eeq 
The 2D restriction and interpolation operators $I_h^{2h}$ and $I_{2h}^h$ are defined as
\beq
I_h^{2h} = I_{h,1D}^{2h} \otimes I_{h,1D}^{2h} \in \mathbb{R}^{\frac{N}{4} \times N}, \qquad 
I_{2h}^h = (I_h^{2h})^T \in \mathbb{R}^{N \times \frac{N}{4}},
\eeq
where $\otimes$ is the Kronecker product. These operators represent a set of first-order restriction and interpolation operators. Note that higher-order intergrid operators may be used for improved accuracy if required. Classical multigrid is based upon the two-grid correction scheme.

\noindent\begin{minipage}{\textwidth}
\vspace{0.2cm}
\noindent \hrulefill \, \textbf{Classical two-grid correction scheme (TG)} \hrulefill \vspace{-0.2cm}
\begin{itemize}
\item[1.] Relax $\nu_1$ times on the equation $A^h x^h = f^h$. \vspace{-0.1cm}
\item[2.] Calculate $r^h = f^h - A^h x^h$ and restrict the residual $r^{2h} = I_h^{2h} r^h$. \vspace{-0.1cm}
\item[3.] Solve the residual equation $A^{2h} e^{2h} = r^{2h}$ for $e^{2h}$ on the coarse grid. \vspace{-0.1cm}
\item[4.] Interpolate the coarse grid error $e^{h} = I_{2h}^h e^{2h}$ to obtain a fine grid error approximation, and correct the initial guess $x^h \leftarrow x^h + e^h$. \vspace{-0.1cm}
\item[5.] Relax $\nu_2$ times on the equation $A^h x^h = f^h$. \vspace{-0.1cm}
\end{itemize}
\vspace{-0.3cm}
\noindent\hrulefill
\vspace{0.2cm}
\end{minipage}

The relaxation in Step 1.~and 5.~applies a basic iterative relaxation scheme (e.g.~weighted Jacobi, Gauss-Seidel, SIRT) to the system. This is commonly referred to as \emph{pre-} and \emph{post-smoothing}. The coarse grid solve in Step 3.~is typically done by recursively embedding the correction scheme, restricting the coarse grid solution to an even coarser grid, etc.,~building a complete multi-level hierarchy. The resulting self-embedded multi-level structure is typically referred to as a multigrid V-cycle. 

\begin{figure}[t]
\begin{center}
\includegraphics[width=12cm]{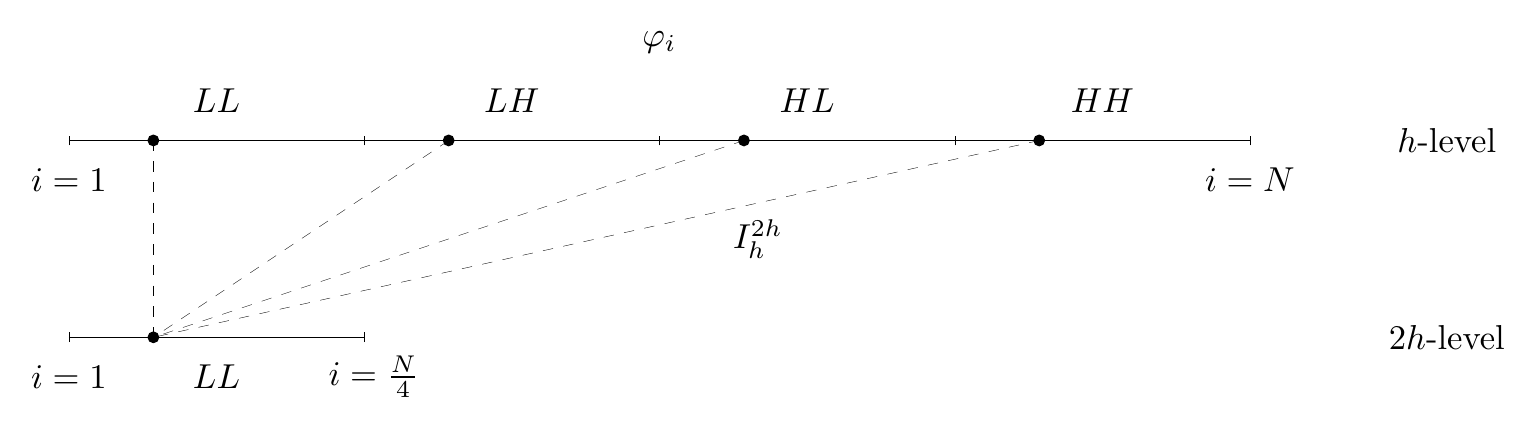}
\caption{Schematic representation of the action of the classical 2D restriction operator $I_h^{2h}$ on the smoothest (LL), semi-oscillatory (LH and HL) and oscillatory (HH) fine grid eigenmodes. Each set of four fine grid harmonics is projected onto a single coarse grid eigenmode in the LL range, see (\ref{eq:coincide}).}
\label{fig:restr1}
\end{center}
\end{figure}

\begin{figure}[t]
\begin{center}
\includegraphics[width=0.48\textwidth]{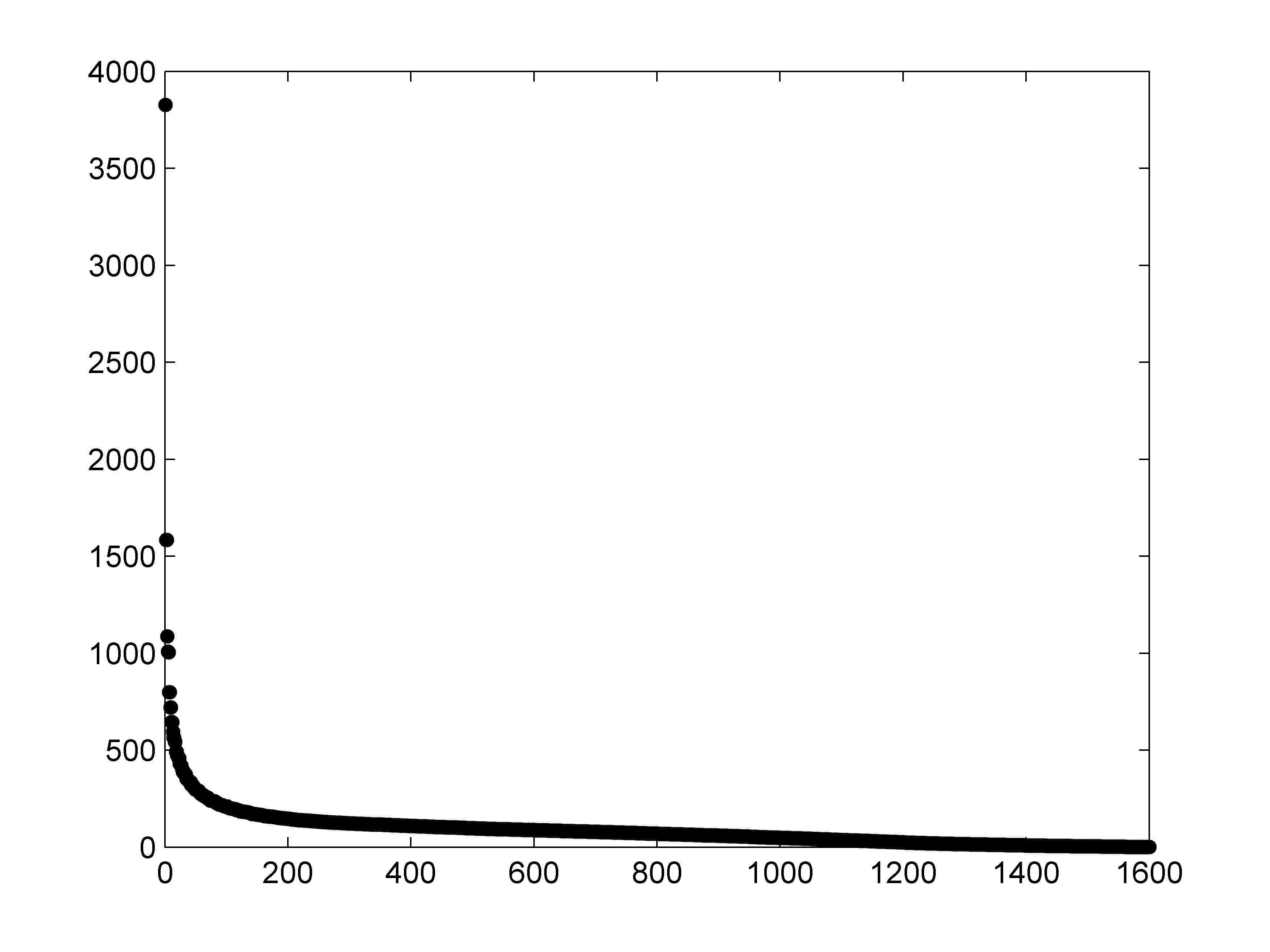}
\includegraphics[width=0.48\textwidth]{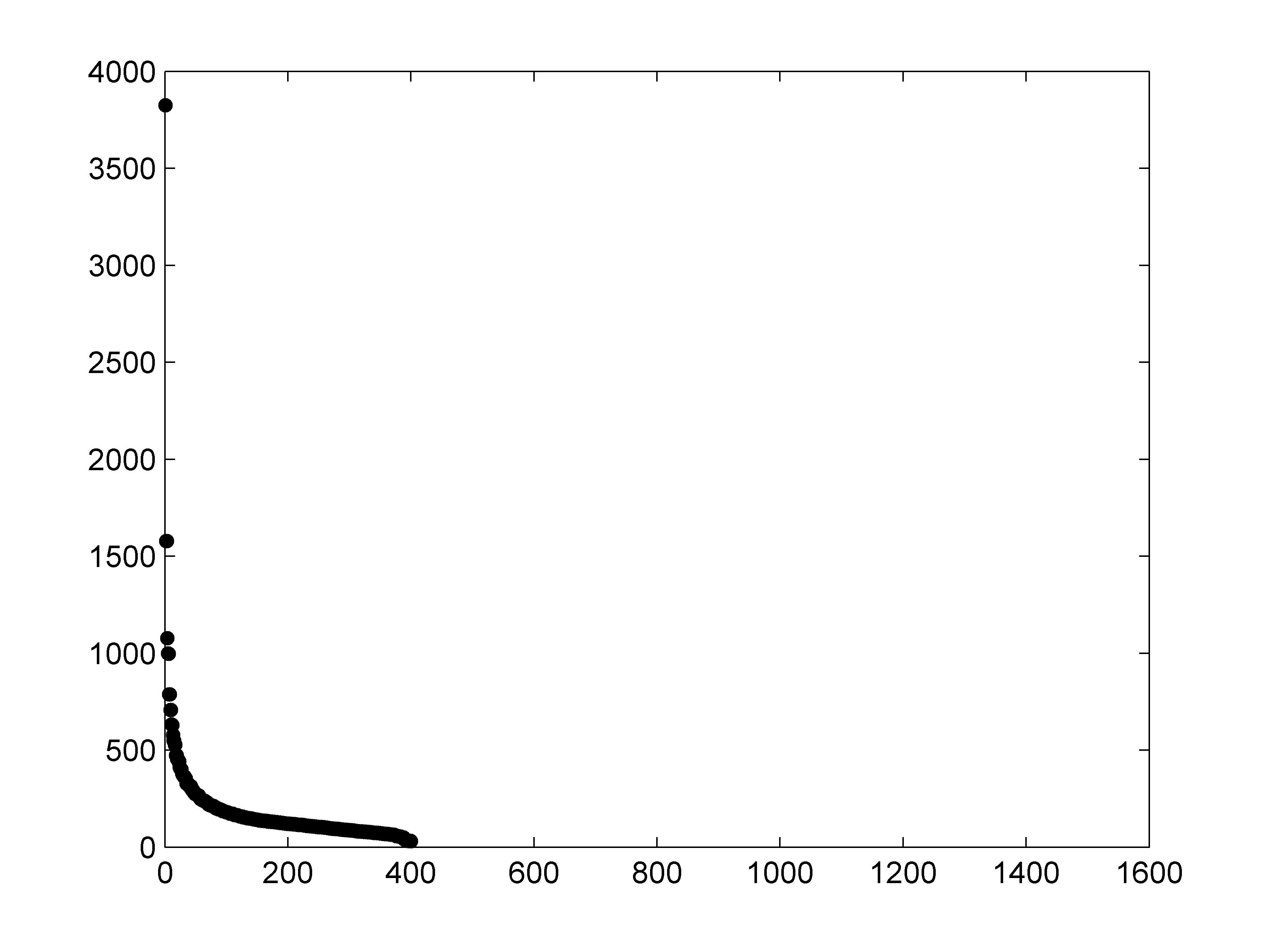}
\caption{Left: spectrum $\lambda^h_i$ $(i=1,\ldots,N)$ of the fine grid operator $A^h = W^T W$, with $N = 40 \times 40$ and $M = 100 \times 40$. Right: spectrum $\lambda^{2h}_i$ $(i=1,\ldots,\frac{N}{4})$ of the coarse grid operator $A^{2h}$, see (\ref{eq:galerkin}).}
\label{fig:eig_val_A}
\end{center}
\end{figure}

For the purpose of analysis, the two-grid correction scheme is commonly considered as an approximation to the entire multigrid cycle \cite{stüben1982multigrid,hackbusch1985multi}. The correction is given by
\beq \label{eq:mg_error}
e^{(k)} = TG \, e^{(k-1)} = S^{\nu_2} (I - I_{2h}^h {A^{2h}}^{-1} I_h^{2h} A^h) S^{\nu_1} \, e^{(k-1)},
\eeq
where the operator $S$ represents the iteration matrix of a basic iterative scheme, and $\nu_1$ and $\nu_2$ are the number of pre- and post-smoothing steps respectively. Sets of so-called \emph{harmonic modes} coincide on $\Omega^{2h}$ due to restriction \cite{briggs2000multigrid}. In 2D, sets of four harmonic modes coincide on the coarse grid \cite{trottenberg2001multigrid,cools2013local},
\beq \label{eq:coincide}
\varphi^{2h}_i = I_{h}^{2h} \varphi^h_{i} = I_{h}^{2h} \varphi^h_{i+\frac{N}{4}} = I_{h}^{2h} \varphi^h_{i+\frac{N}{2}} = I_{h}^{2h} \varphi^h_{i+3\frac{N}{4}}, \quad i = 1,\ldots,\frac{N}{4}.
\eeq
Note that the eigenmodes in (\ref{eq:coincide}) are ordered in function of the sets of harmonics for notational convenience. This ordering may differ slightly from the sorted eigenvalue-based ordering introduced in Section \ref{sec:classical_art}. The concept of coinciding modes is visualized on Figure \ref{fig:restr1}.
A well-known result from multigrid theory, see \cite{briggs2000multigrid,hackbusch1985multi}, states that for $k \in \{0,1,2,3\}$
\begin{equation}
(I - I_{2h}^h {A^{2h}}^{-1} I_h^{2h} A^h) \, \varphi^h_{i+k\frac{N}{4}} \approx 
\left(1 - \frac{\lambda^h_{i+k\frac{N}{4}}}{\lambda^{2h}_i} \right) \varphi^h_{i+k\frac{N}{4}}, \quad i = 1,\ldots,\frac{N}{4} , \label{eq:briggs_smooth}
\end{equation}
where $\lambda^h_i$ denotes the eigenvalue of $A^h$ corresponding to the $i$-th eigenmode $\varphi^h_i$.
Under the assumption that $\lambda^h_i \approx \lambda^{2h}_i$ for the smooth eigenmodes, which generally holds (see Figure \ref{fig:eig_val_A}), Eqn.~(\ref{eq:briggs_smooth}) implies that the smoothest modes are approximately mapped onto zero by the two-grid operator. For oscillatory modes, however, it holds that $\lambda^h_{i+k\frac{N}{4}} \ll \lambda^{2h}_i$ $(k \in \{1,2,3\})$, hence oscillatory modes are left unchanged. This in principle opposes the action of the smoother $S$, which is assumed to damp oscillatory modes while leaving smooth modes unaffected. The complementary action of smoother and two-grid correction is crucial for the effectiveness of the multigrid scheme. 

Although multigrid methods allow for fast and scalable solving and/or preconditioning of a wide range of PDE-type problems, this is generally not true for tomographic reconstruction. As shown by R\"ude et al.~in their work on AMG for ART acceleration \cite{kostler2006towards}, multigrid preconditioning does not significantly reduce the number of Krylov iterations. This ineffectiveness is evident from the discussion above, since basic iterative relaxation schemes do not possess the smoothing property for algebraic tomographic reconstruction problems. Hence, the classical multigrid scheme induces a damping of the smooth modes by \emph{both} the smoother and correction operator while oscillatory modes remain present in the error, causing reduced convergence and performance.

\section{Wavelet-based multigrid preconditioner}

In this section we introduce a new multi-level preconditioner for Krylov methods, which is specifically tailored to tomographic reconstruction problems. Contrarily to standard multigrid preconditioners, the proposed scheme does not rely on the smoothing properties of a basic relaxation scheme. Instead, we aim at constructing a multi-level scheme in which the damping of all eigenmodes is incorporated within the correction scheme itself, and the intergrid operators are adapted to this purpose. We introduce different wavelet-based operators for the various spectral regions, which allows elimination of all error components by consecutive coarse grid projection.

\subsection{Definition and notation}

We define a more advanced multigrid correction scheme, inspired by the theory of scaling- and wavelet-functions.

\subsubsection{Intergrid operators.} 
\begin{figure}[t]
\begin{center}
\includegraphics[width=12cm]{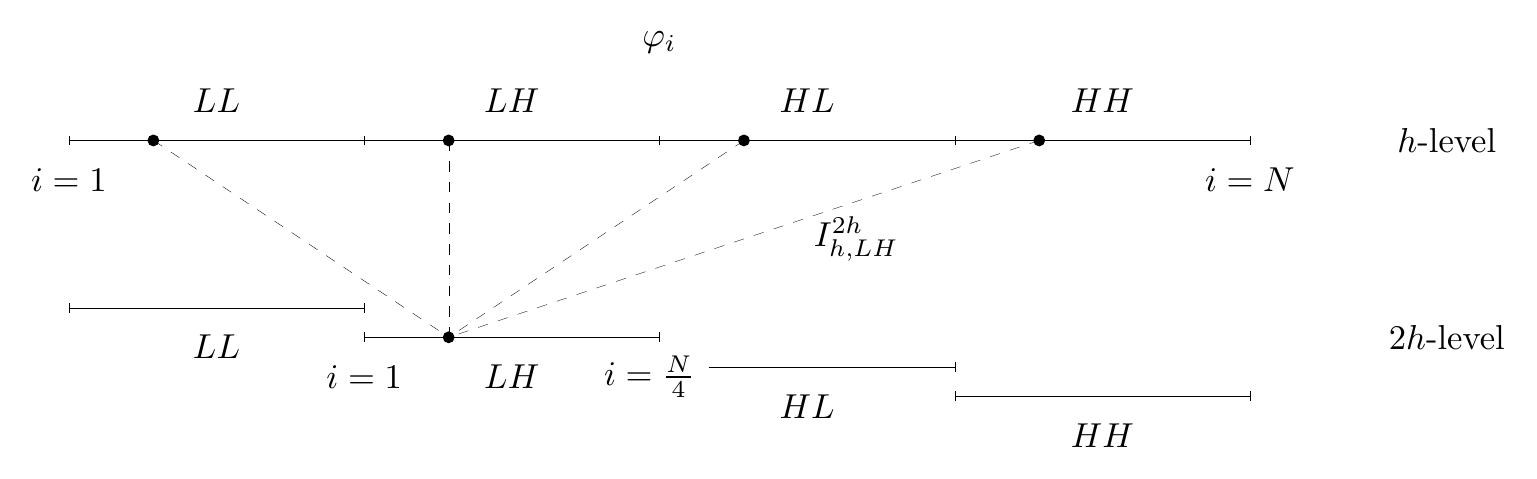}
\caption{Schematic representation of the action of the 2D wavelet-based restriction operator $I_{h,LH}^{2h}$ on the smoothest (LL), semi-oscillatory (LH and HL) and oscillatory (HH) fine grid eigenmodes. Each set of four fine grid harmonics is projected onto a single coarse grid mode in the LH range, see (\ref{eq:coincide2}).}
\label{fig:restr2}
\end{center}
\end{figure}

Let the basic 1D intergrid operators $I_{h,1D}^{2h}$ and $I_{2h,1D}^h$ be defined by (\ref{eq:intops}). The rows of $I_{h,1D}^{2h}$ are commonly referred to as discrete Haar \emph{scaling functions} in the wavelet literature \cite{haar1910theorie}. Additionally, we define the wavelet operators $J_{h,1D}^{2h}$ and $J_{h,1D}^{2h}$ based upon the Haar \emph{wavelet functions} corresponding to these scaling functions as
\beq \label{eq:intops_wav}
J_{h,1D}^{2h} = \frac{1}{\sqrt{2}}
\begin{bmatrix} 
1 & \hspace{-0.2cm} -1 & & & & & &\\
& & 1 & \hspace{-0.2cm} -1 & & & & \\
& & & & \cdot &\cdot & &\\
& & & & & & 1 & \hspace{-0.2cm} -1 \\
\end{bmatrix} \in \mathbb{R}^{\frac{n}{2} \times n},
~~
J_{2h,1D}^{h} = (J_{h,1D}^{2h})^T \in \mathbb{R}^{n \times \frac{n}{2}}.
\eeq 
Higher order scaling and wavelet functions, e.g.~D4 Daubechies functions \cite{daubechies1988orthonormal}, may be used to replace the Haar functions to obtain a higher precision, but are computationally more expensive. Hence, we restrict ourselves to the simplest class of Haar-type scaling and wavelet functions. Using these one-dimensional scaling and wavelet operators, four sets of 2D intergrid operators are defined:
\begin{align}
I_{h,LL}^{2h} &= I_{h,1D}^{2h} \otimes I_{h,1D}^{2h} \in \mathbb{R}^{\frac{N}{4} \times N}, \qquad 
I_{2h,LL}^h = (I_{h,LL}^{2h})^T \in \mathbb{R}^{N \times \frac{N}{4}}, \label{eq:1of4} \\
I_{h,LH}^{2h} &= I_{h,1D}^{2h} \otimes J_{h,1D}^{2h} \in \mathbb{R}^{\frac{N}{4} \times N}, \qquad 
I_{2h,LH}^h = (I_{h,LH}^{2h})^T \in \mathbb{R}^{N \times \frac{N}{4}}, \\
I_{h,HL}^{2h} &= J_{h,1D}^{2h} \otimes I_{h,1D}^{2h} \in \mathbb{R}^{\frac{N}{4} \times N}, \qquad 
I_{2h,HL}^h = (I_{h,HL}^{2h})^T \in \mathbb{R}^{N \times \frac{N}{4}}, \\
I_{h,HH}^{2h} &= J_{h,1D}^{2h} \otimes J_{h,1D}^{2h} \in \mathbb{R}^{\frac{N}{4} \times N}, \qquad 
I_{2h,HH}^h = (I_{h,HH}^{2h})^T \in \mathbb{R}^{N \times \frac{N}{4}}. \label{eq:4of4}
\end{align}
As illustrated by Figure \ref{fig:restr2}, these operators project the eigenspace of $A^h$ onto four disjunct coarse grid subspaces, designated by the subscript indices $LL$, $LH$, $HL$ and $HH$. The first restriction operator $I_{h,LL}^{2h}$ is the standard restriction, which maps sets of four harmonic eigenmodes onto a single smooth coarse grid representative, see Eqn.~(\ref{eq:coincide}) and Figure \ref{fig:restr1}. The operators $I_{h,LH}^{2h}$ and $I_{h,HL}^{2h}$ map the harmonics onto a coarse grid representative which is slowly varying in the $x$-direction but oscillatory in the $y$-direction, or vice versa. The action of $I_{h,LH}^{2h}$ is illustrated on Figure \ref{fig:restr2}. Finally, $I_{h,HH}^{2h}$ maps the harmonic eigenmodes onto a subset of coarse grid representatives which are highly oscillatory in all directions. Formally we write, in analogy to (\ref{eq:coincide}),
\beq \label{eq:coincide2}
\varphi^{2h}_{i,id} = I_{h,id}^{2h} \, \varphi^h_{i} = I_{h,id}^{2h} \, \varphi^h_{i+\frac{N}{4}} = I_{h,id}^{2h} \, \varphi^h_{i+\frac{N}{2}} = I_{h,id}^{2h} \, \varphi^h_{i+3\frac{N}{4}}, \quad i = 1,\ldots,\frac{N}{4},
\eeq
for $id \in \{LL,LH,HL,HH\}$, where $\varphi^{2h}_{i,id}$ is the $i$-th coarse grid eigenmode in the respective range of the restriction operator $I_{h,id}^{2h}$ .

\subsubsection{Wavelet-based two-grid scheme.} \label{sec:wavelet_twogrid}

We now define a wavelet-based two-grid scheme that performs four (consecutive) two-grid correction steps, using the sets of intergrid operators defined above.

\noindent\begin{minipage}{\textwidth}
\vspace{0.2cm}
\noindent \hrulefill \, \textbf{Wavelet-based two-grid correction scheme (WTG)} \hrulefill \vspace{0.1cm} \newline
for $id \in \{LL,LH,HL,HH\}$ \vspace{-0.1cm}
\begin{itemize}
\item[1.] Calculate $r^h = b^h - A^h x^h$ and project to a coarse grid $r_{id}^{2h} = I_{h,id}^{2h} r^h$. \vspace{-0.1cm}
\item[2.] Solve the residual equation $A_{id}^{2h} e_{id}^{2h} = r_{id}^{2h}$ for $e_{id}^{2h}$ on the coarse grid. \vspace{-0.1cm}
\item[3.] `Interpolate' the coarse grid error $e_{id}^{h} = I_{2h,id}^h e_{id}^{2h}$ to obtain a fine grid error approximation, and correct the initial guess $x^h \leftarrow x^h + e_{id}^h$. \vspace{-0.1cm}
\end{itemize}
\vspace{-0.3cm}
\noindent\hrulefill
\vspace{0.3cm}
\end{minipage}

\noindent The coarse grid operators $A_{id}^{2h}$ are defined using Galerkin products
\beq
A_{id}^{2h} = I_{h,id}^{2h} \, A^h \, I_{2h,id}^h, \quad id \in \{LL,LH,HL,HH\}.
\eeq
In terms of error propagation, the wavelet two-grid scheme can be written as
\begin{align} \label{eq:wmg}
e^{(k)} = & \, WTG \, e^{(k-1)} \notag \\
= & \, (I - I_{2h,HH}^h (A_{HH}^{2h})^{-1} I_{h,HH}^{2h} A^h) (I - I_{2h,HL}^h (A_{HL}^{2h})^{-1} I_{h,HL}^{2h} A^h) \\ 
  & \, (I - I_{2h,LH}^h (A_{LH}^{2h})^{-1} I_{h,LH}^{2h} A^h) (I - I_{2h,LL}^h (A_{LL}^{2h})^{-1} I_{h,LL}^{2h} A^h) \, e^{(k-1)} \notag
\end{align}
This correction scheme solves the error equation by consecutive projection onto the four coarse grid subspaces consisting of smooth ($LL$), semi-oscillatory ($LH$ and $HL$) and oscillatory ($HH$) modes. This leads to an elimination of smooth, semi-oscillatory and highly oscillatory error modes respectively, thus resolving the entire error spectrum in every iteration.
Note that no relaxation is incorporated in the WMG scheme, since standard relaxation methods perform poorly for tomographic reconstruction problems. Instead, damping of the oscillatory eigenmodes is accomplished in a natural way by projection onto the oscillatory coarse grid subspaces $LH$, $HL$ and $HH$. Analogously to (\ref{eq:briggs_smooth}), the projection onto the oscillatory subspaces satisfies the following relations for $k \in \{0,1,2,3\}$ and $i = 1,\ldots,\frac{N}{4}$:
\begin{equation}
(I - I_{2h,id}^h (A_{id}^{2h})^{-1} I_{h,id}^{2h} A^h) \, \varphi^h_{i+k\frac{N}{4}} \approx 
\left(1 - \frac{\lambda^h_{i+k\frac{N}{4}}}{\lambda^{2h}_{i,id}} \right) \varphi^h_{i+k\frac{N}{4}}. \label{eq:briggs_oscil_H}
\end{equation}
Observing that for oscillatory eigenmodes ($k\in\{1,2,3\}$) we have $\lambda^h_{i+k\frac{N}{4}} \approx \lambda^{2h}_{i,id}$ ($id \in \{LH,HL,HH\}$), Eqn.~(\ref{eq:briggs_oscil_H}) implies oscillatory modes are eliminated by the projection onto the oscillatory coarse grid subspaces, while smooth modes are left unchanged since $\lambda^h_i \ll \lambda^{2h}_{i,LL}$. This indicates that the projections onto the oscillatory coarse grid subspaces $LH$, $HL$ and $HH$ indeed eliminate the oscillatory eigencomponents from the error.


Note that the above definition of the WTG scheme implies a multiplicative multi-level formulation, as the residual is recalculated in every step using the corrected guess $x^h$. Alternatively, the residual calculation may be performed outside the loop, leading to an additive variant of the WTG scheme, which is ideally suited for multi-core parallelization yet generally features reduced convergence speed compared to the multiplicative variant. In this work we have opted for a hybrid approach, where one residual calculation is performed after the elimination of the smooth components by the $LL$ projection step. Since the initial error guess on every level is zero, this yields a computational cost of one SpMV per projection, analogous to standard multigrid, while retaining optimal stability. Additionally, our approach allows for a parallel execution of the three latter projections ($LH$, $HL$ and $HH$). For the numerical experiments in this work, however, we have alternatively chosen to parallelize the (dense) coarse grid exact solves using BLAS3 routines, instead of parallellizing over the coarse grid projection spaces.

\subsubsection{Wavelet-based multigrid (WMG).} A WMG V-cycle consists of a recursive embedding of the WTG two-grid scheme. Denoting the total number of levels in the WMG hierarchy by $\ell$ and assuming that $A^h \in \mathbb{R}^{N \times N}$ as above, the original fine grid problem $A^h x^h = f^h$ is split up into a collection of $4^{\ell-1}$ subproblems of size $\frac{N}{4^{\ell-1}} \times \frac{N}{4^{\ell-1}}$ on the coarsest grid, which are computationally much cheaper to solve. Note that the coarse grid operators $A^{2h}_{id}$ feature the same sparsity structure over all levels (including identical nonzeros-zeros ratio) due to the choice of the projector basis functions (\ref{eq:1of4})-(\ref{eq:4of4}). 

The advantage of WMG as a Krylov preconditioner over other preconditioning techniques like e.g.~incomplete Cholesky factorization \cite{polydorides2002krylov} is that instead of directly solving the large-scale fine grid system, the problem is reorganized towards solving a collection of small subproblems, which are more amenable to direct solution.
Note that this key idea of the WMG method resembles the Hierarchical Basis Multigrid Method (HBMM) \cite{bank1988hierarchical,bank1996hierarchical}. For both techniques, the division of the large scale system into a collection of smaller-basis subproblems allows for a direct solution of the coarse grid problems. Moreover, this collection of subproblems is particularly suited for multi-core parallelization (see the discussion on parallellization in Section \ref{sec:wavelet_twogrid} above). 

The WMG method can in principle be applied as a stand-alone solver to the system (\ref{eq:sys})-(\ref{eq:sys_ne}). In this work we however opt to use the WMG scheme as a Krylov preconditioner, since the embedding in a governing Krylov solver generally leads to a faster and more robust solution scheme. Note that the preconditioner is approximately inverted using only one WMG cycle, as is common practice in the MG literature. In the next sections the efficiency of the WMG-preconditioned Krylov method is analyzed. 
We again stress that the matrix $A^h = W^T W$ and its coarse grid representations $A_{id}^{2h} = I_{h,id}^{2h} \, W^T \, W \, I_{2h,id}^h$ are never computed explicitly, as this would result in a dense matrix operator with a large memory footprint. Instead, the `tall-and-skinny' sparse matrices $(W \, I_{2h,id}^h) \in \mathbb{R}^{M \times N/4}$ are computed and stored on each level in a preliminary setup step, and can hence be applied as an SpMV operation at any point in the algorithm. 

Finally, we briefly comment on the computational complexity of one WMG cycle. Assuming the fine level SpMV operation features a cost of $\mathcal{O}(N)$, we note that each coarse level SpMV operation is four times cheaper. However, since the number of coarse grid subproblems increases by a factor 4 on each coarser level, the per-level cost of the WMG scheme is a constant $\mathcal{O}(N)$ operations. Hence, since the number of levels is proportionate to $\mathcal{O}(\log N)$, the total computational cost of one WMG cycle is $\mathcal{O}(N \log N)$.

\subsection{Spectral properties of the WMG-Krylov method}

\begin{figure}[t]
\begin{center}
\includegraphics[width=0.32\textwidth]{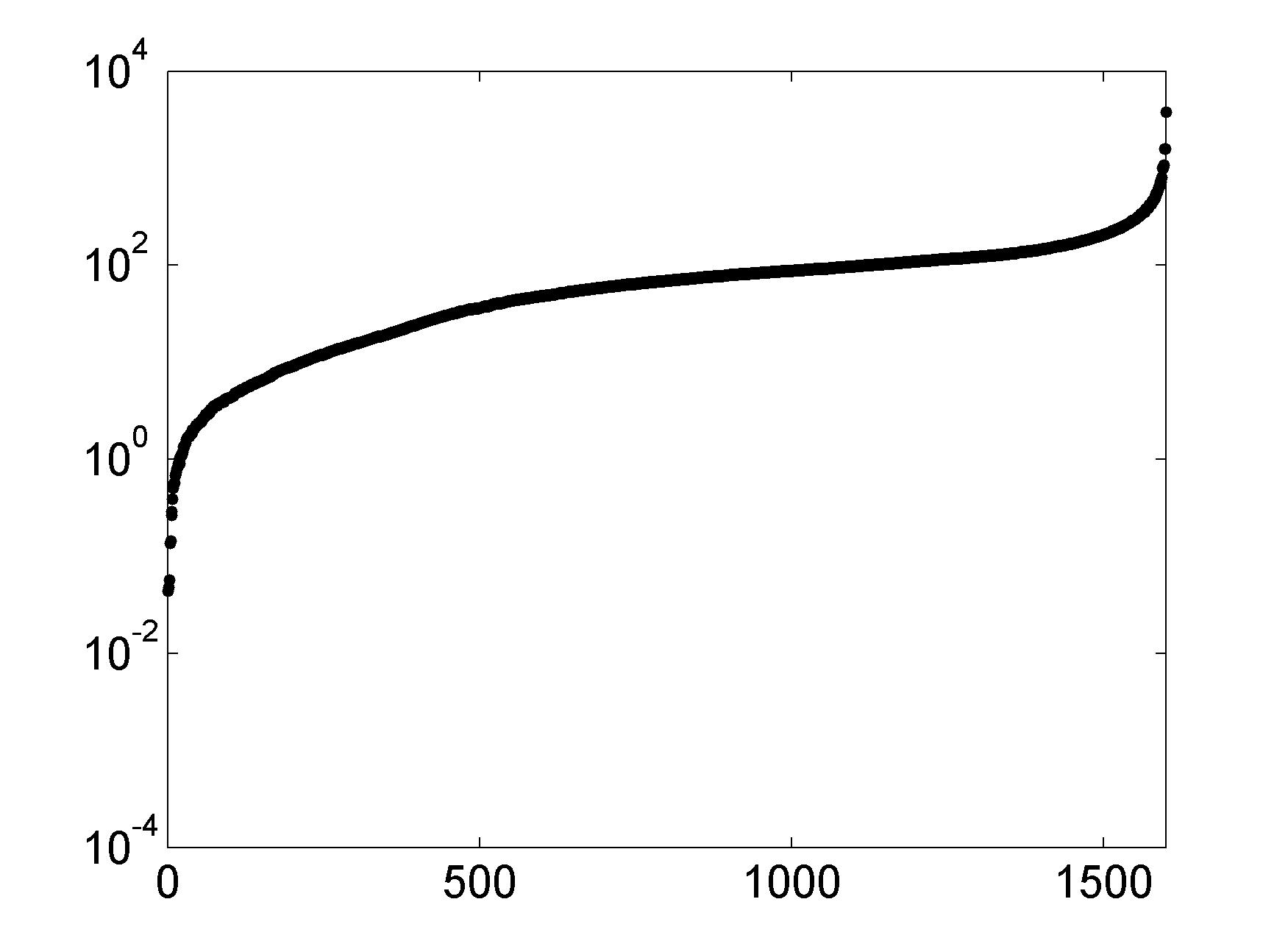}
\includegraphics[width=0.32\textwidth]{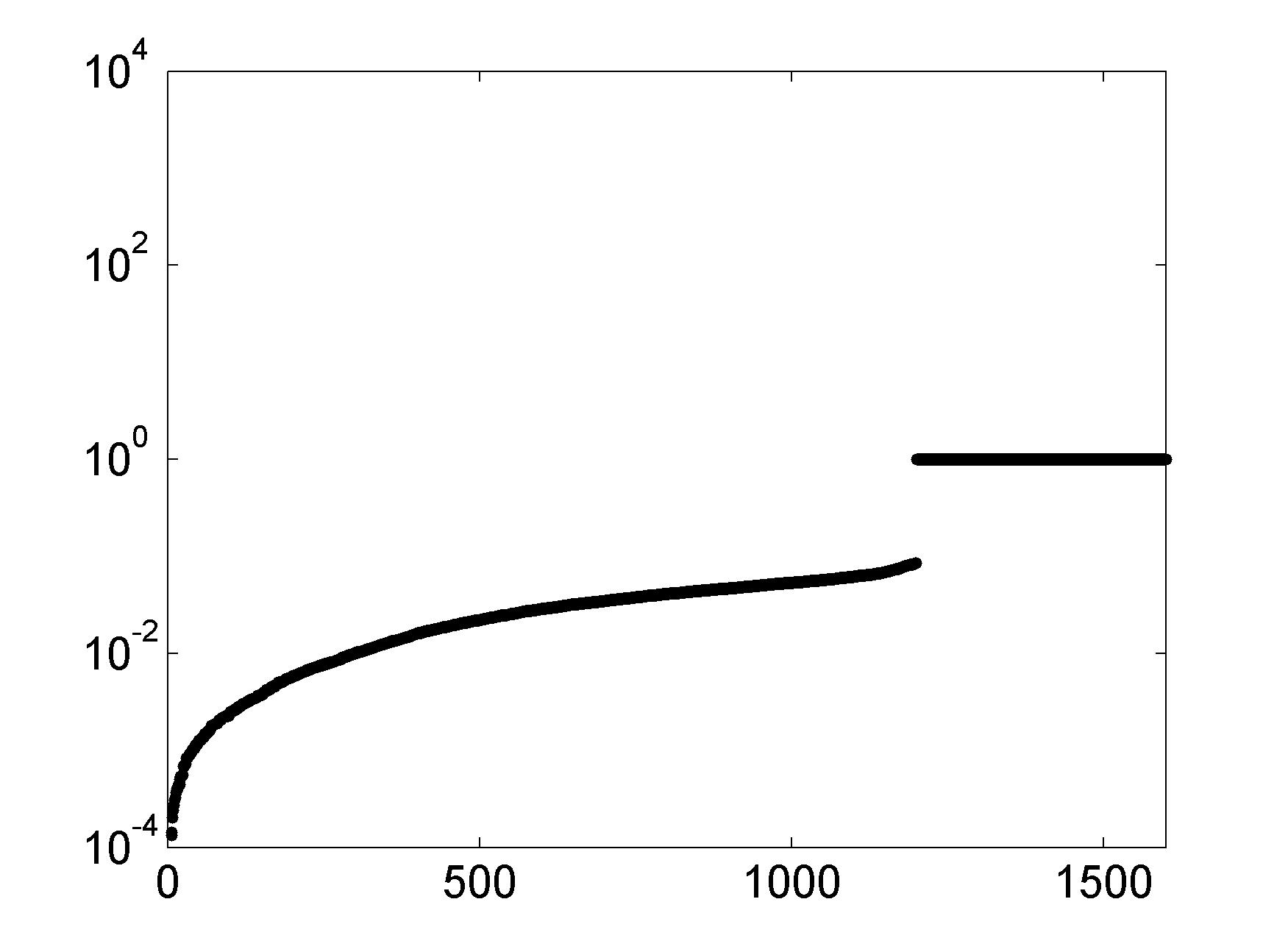}
\includegraphics[width=0.32\textwidth]{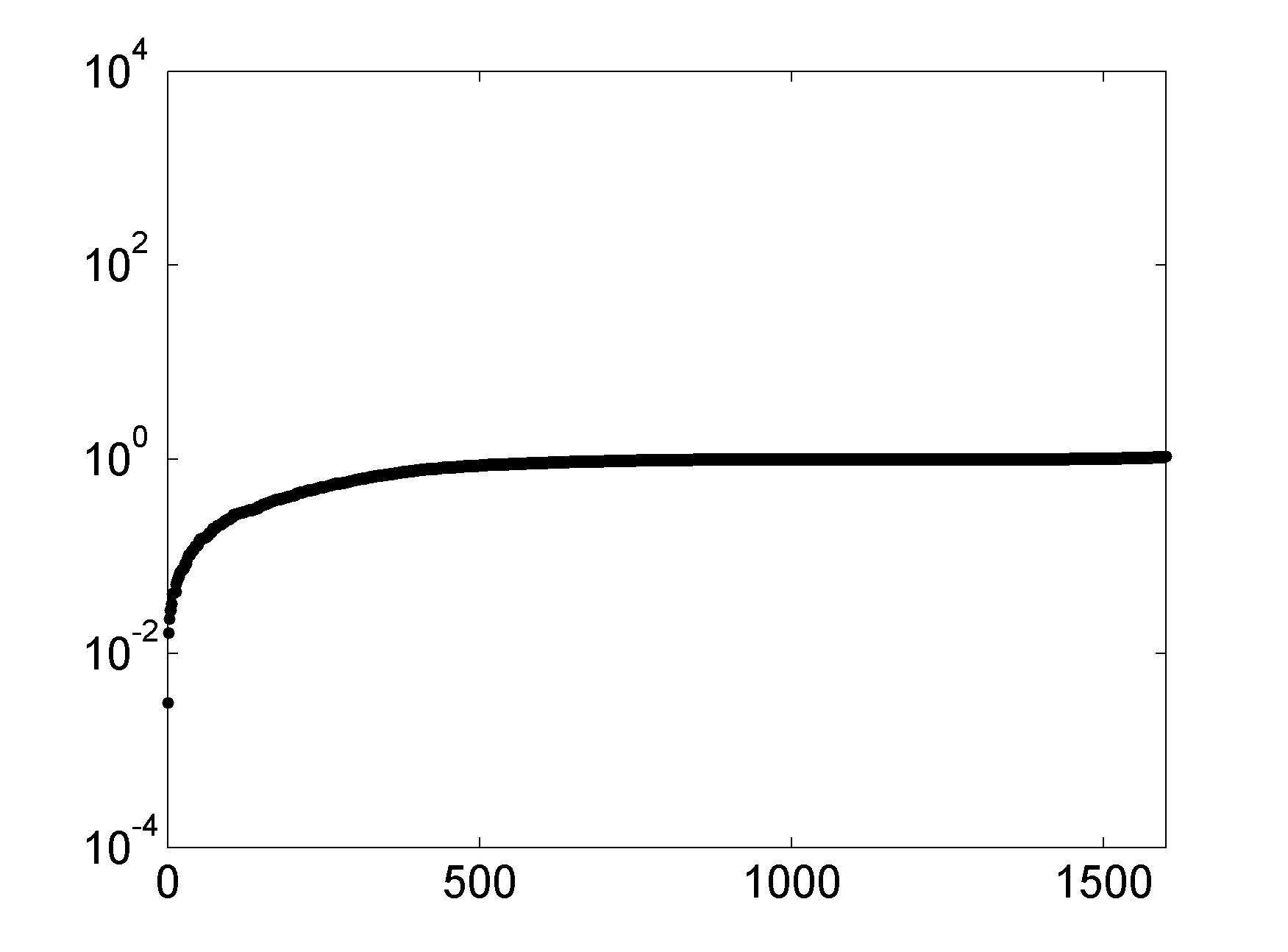}
(a) \hspace{0.28\textwidth} (b) \hspace{0.28\textwidth} (c)
\caption{Spectra of the $N = 40 \times 40$ by $M = 100 \times 40$ iteration matrices of (a) the unpreconditioned Krylov method (\ref{eq:spec_vanilla}) with corresponding condition number $\kappa = 8.68e4$, (b) the TG-Krylov method (\ref{eq:spec_classic}) with $\kappa = 4.52e4$, and (c) the WTG-Krylov method (\ref{eq:spec_wavelet}) with $\kappa = 3.42e2$. 
Vertical axis in log scale.} 
\label{fig:3specs}
\end{center}
\end{figure}

To obtain more insight in the potential convergence improvement when using the WMG scheme as a Krylov preconditioner (WMG-Krylov), we analyze the spectral properties of the WTG-Krylov method. The analysis presented here is based on the analysis of multigrid-preconditioned Krylov methods performed in \cite{cools2013local} and \cite{wienands2000fourier}. The spectrum of the original system matrix $W^T W$,
\beq \label{eq:spec_vanilla}
\sigma(W^T W),
\eeq
is compared to the spectrum of the WTG-preconditioned Krylov iteration matrix (right preconditioning), given by
\beq \label{eq:spec_wavelet}
\sigma(I-(W^T W) \, WTG \,  (W^T W)^{-1}),
\eeq
see \cite{wienands2000fourier}, page 590, Eqn.~(21), where the operator $WTG$ is defined by (\ref{eq:wmg}). For completeness of comparison, we also consider the spectrum of the standard TG-preconditioned Krylov iteration matrix
\beq \label{eq:spec_classic}
\sigma(I-(W^T W) \, TG \, (W^T W)^{-1}),
\eeq
with $TG$ defined by (\ref{eq:mg_error}).
Figure \ref{fig:3specs} shows the spectrum of the unpreconditioned operator $W^T W$ (\ref{eq:spec_vanilla}), the operator preconditioned by classical $TG(1,1)$ two-grid with one pre- and post-smoothing SIRT iteration (\ref{eq:spec_classic}) and the operator preconditioned by the WTG correction scheme (\ref{eq:spec_wavelet}) for a volume size of $N = 40 \times 40$ with 100 equiangular parallel beam projections (over $180^\circ$) of 40 rays each. 
Note how only the WMG-preconditioned spectrum is efficiently clustered around one, showing it to be much more amenable to Krylov solution than the unpreconditioned or MG-preconditioned systems. Indeed, the condition number
\beq \label{eq:condnr}
\kappa(A) = \frac{\max_i \, \left|\lambda_i(A)\right|}{\min_i \, \left|\lambda_i(A)\right|},  \quad  (i = 1,\ldots,N),
\eeq
is reduced significantly by the WMG preconditioner, lying in the order of magnitude $\mathcal{O}(10^{2})$ compared to $\mathcal{O}(10^{4})$ for the unpreconditioned and MG-Krylov method. Hence, the WMG-Krylov method is expected to converge significantly faster than the unpreconditioned Krylov method.

\section{Numerical results}

In this section a variety of 2D benchmark problems is presented to compare the performance of the proposed WMG-Krylov solver to the SIRT iterative reconstruction technique and the unpreconditioned Krylov solver. We first consider a simple test case with noiseless data to validate the effectiveness of the WMG preconditioner. 
Subsequently, a more realistic benchmark problem with the addition of an artificial random noise component to the data is analyzed. 

\subsection{Shepp-Logan type model problem with non-noisy data}
\begin{figure}[t]
\begin{center}
\includegraphics[width=1.0\textwidth]{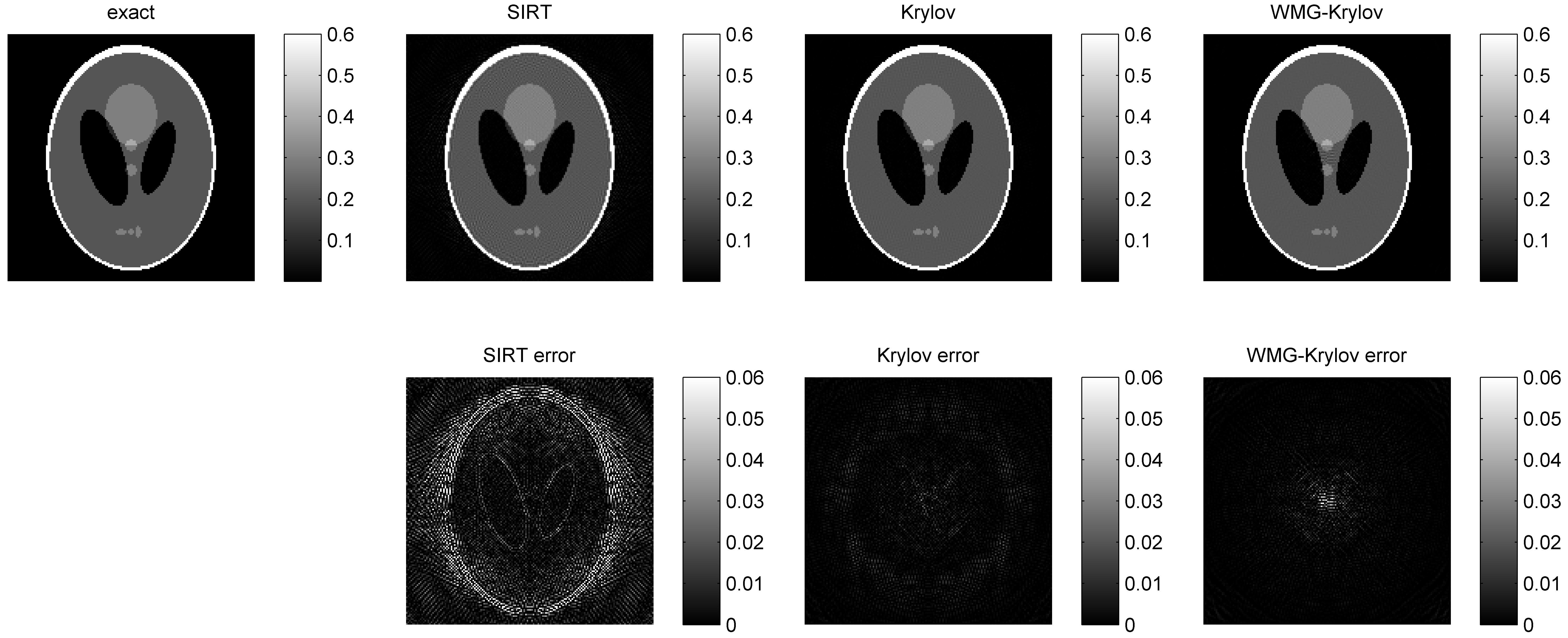}
\caption{Shepp-Logan type model problem with $N = 160 \times 160$ and $M = 400 \times 160$ (no noise). Displayed are (f.l.t.r.)~the exact solution $x_{ex}$ and numerical solutions computed by 1000 SIRT iterations, 300 BiCGStab iterations and 50 WMG-preconditioned BiCGStab iterations respectively. Specifications: see Table \ref{tab:results}.}
\label{fig:sol}
\end{center}
\end{figure}

For the first test case we restrict ourselves to a model problem with non-noisy data, aiming at an exact reconstruction of the image. Consider a Shepp-Logan type model problem consisting of a 160-by-160 pixel image. The exact solution to this model problem is denoted by $x_{ex}$ (see Figure \ref{fig:sol}, top left) and is known explicitly, allowing to calculate the error norm $\|e_k\| = \|x_k - x_{ex}\|$ after every Krylov iteration. The Shepp-Logan phantom object of interest is projected under 400 equiangular projection angles (equally distributed over $180^\circ$) using 160 rays per angle, inducing a complete dataset since the object of interest is contained within the inner circle of the image. 

\begin{figure}[t]
\begin{center}
\includegraphics[width=0.6\textwidth]{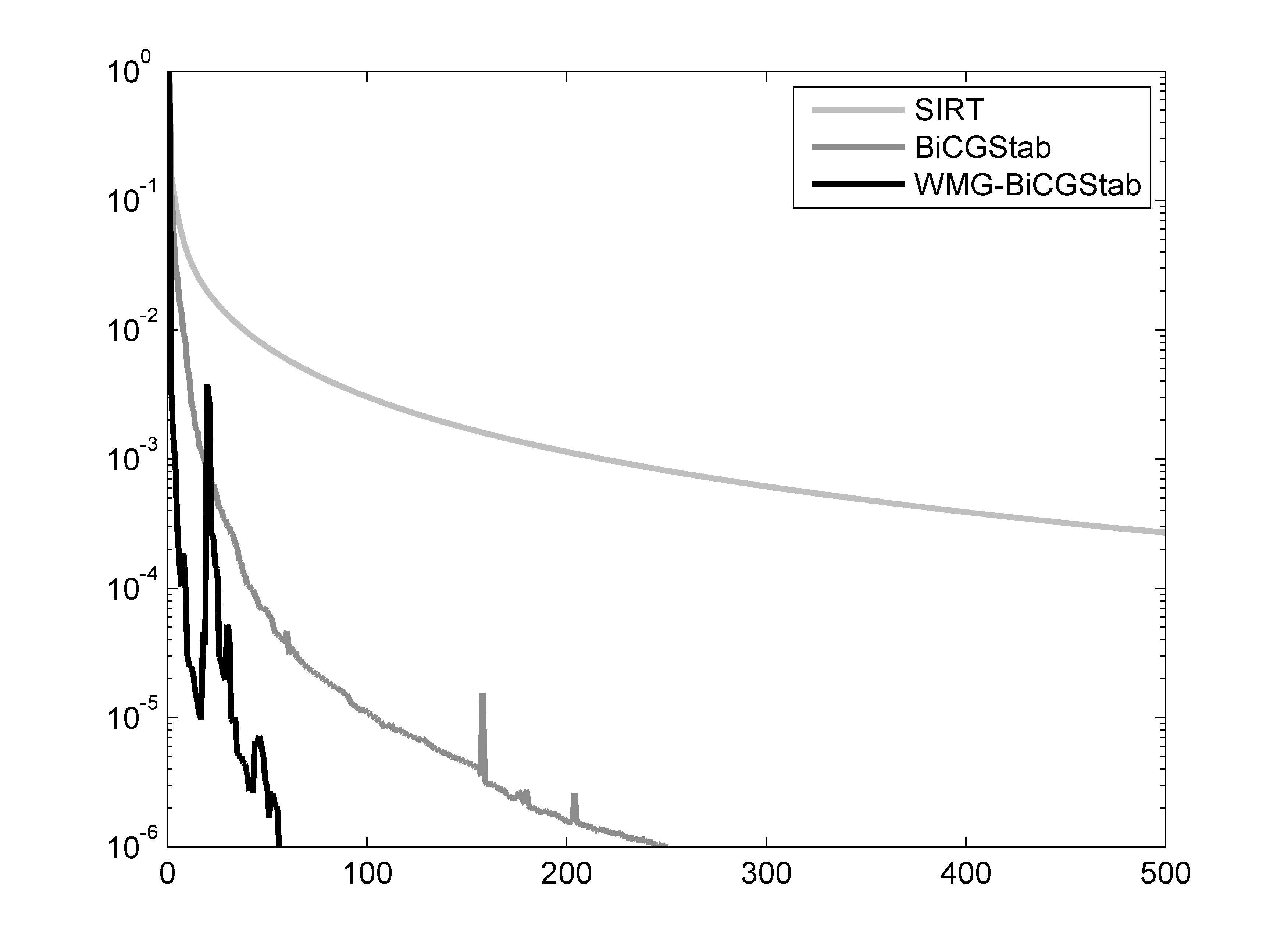}
\caption{Shepp-Logan type model problem with $N = 160 \times 160$ and $M = 400 \times 160$ (no noise). Displayed are the relative residual $L_2$ norms $\|r_k\|_2/\|r_0\|_2$ in function of the number of iterations for SIRT (light gray), BiCGStab (dark gray) and WMG-BiCGStab (black). Vertical axis in log scale.}
\label{fig:relres}
\end{center}
\end{figure}

We compare the convergence of the WMG-BiCGStab solver to the SIRT and unpreconditioned BiCGStab algorithms. A three-level WMG scheme ($l = 3$) is used as a preconditioner, and the collection of 40-by-40 pixel coarse level subproblems is solved exactly using Cholesky factorization. Figure \ref{fig:relres} shows the residual history for all methods. Shown is the scaled $L_2$ norm of the residual $\|r_k\|_2/\|r_0\|_2$ in every iteration. Convergence of the SIRT scheme tends to stagnate rapidly, and the method is clearly outperformed by the Krylov solvers. Furthermore, the WMG-Krylov method displays a significant speedup compared to the unpreconditioned Krylov solver. The WMG preconditioner considerably improves Krylov convergence speed, as was predicted by the spectral analysis in Section 3.2.

The resulting reconstructions are shown in Figure \ref{fig:sol}. Note that the number of iterations was capped at 1000 and is subject to the relative error criterion $\|e_k\|_2/\|e_0\|_2 < \text{tol}$, where we have chosen $\text{tol} = 0.02$. Details on the reconstruction can be found in Table \ref{tab:results}. The SIRT solution displays a large number of small-scale artifacts resulting in an $L_2$ error of approximately 10\% after 1000 iterations. The unpreconditioned Krylov method generates less artifacts, yielding an $L_2$ error of 1.7\% after 300 BiCGStab iterations. Although the WMG-Krylov solution displays some artifacts near the center of the image, the overall reconstruction is very good, displaying an $L_2$ error of 1.5\% after only 50 iterations. 

\begin{table}[t]
\centering
\begin{tabular}{|c| c c c|}
\hline
  & iterations & CPU time & $L_2$ error \\
\hline
  {\footnotesize SIRT} 				 & 1000+ & 80.6 s. & 0.1015 \\
  {\footnotesize BiCGStab} 		 & 300   & 25.1 s. & 0.0166 \\
  {\footnotesize WMG-BiCGStab} & 50    & 17.4 s. & 0.0152 \\
\hline
\end{tabular}
\caption{Shepp-Logan type model problem with $N = 160 \times 160$ and $M = 400 \times 160$ (no noise). Displayed are the number of iterations required to reach the error criterion (tol = 2\%), elapsed CPU time (in seconds), and relative $L_2$ error $\|e_k\|_2/\|e_0\|_2$ for different solution methods.}
\label{tab:results}
\end{table}

Note that Krylov methods are particularly good at reconstructing the sharp, high-contrast edges of the image, as opposed to the SIRT method which generally tends to smear out sharp edges through consecutive iterations. This is reflected in the relative $L_{\infty}$ norm of the error $\|e_k\|_{\infty}/\|e_0\|_{\infty}$, which is $3.17\%$ and $6.69\%$ for the Krylov and WMG-Krylov methods respectively, compared to $20.1\%$ for SIRT. 
The total computational cost of the WMG-Krylov method is significantly lower than the cost of the SIRT method, see CPU timings.\footnote{System specifications: Intel Core i7-2720QM 2.20GHz CPU, 6MB Cache, 8GB RAM.} Although the per-iteration computational cost of the WMG-Krylov method is higher than the cost of a single unpreconditioned Krylov iteration, the reduced number of Krylov iterations pays off in terms of total computational cost.

The small artifact in the center of the WMG-Krylov reconstruction is due to the action of the WMG scheme, which resolves \emph{all} error eigenmodes, see Section 3. This causes the backprojection to display a small and natural accumulation of high-oscillatory artifacts near the center of the image, which is reflected in the relative $L_{\infty}$ norm of the WMG-Krylov error. This artifact is however naturally resolved by the incorporation of a minor regularization term in the system, as shown in Section 4.2.

\subsection{Regularization}

\begin{figure}[t]
\begin{center}
\includegraphics[width=0.49\textwidth]{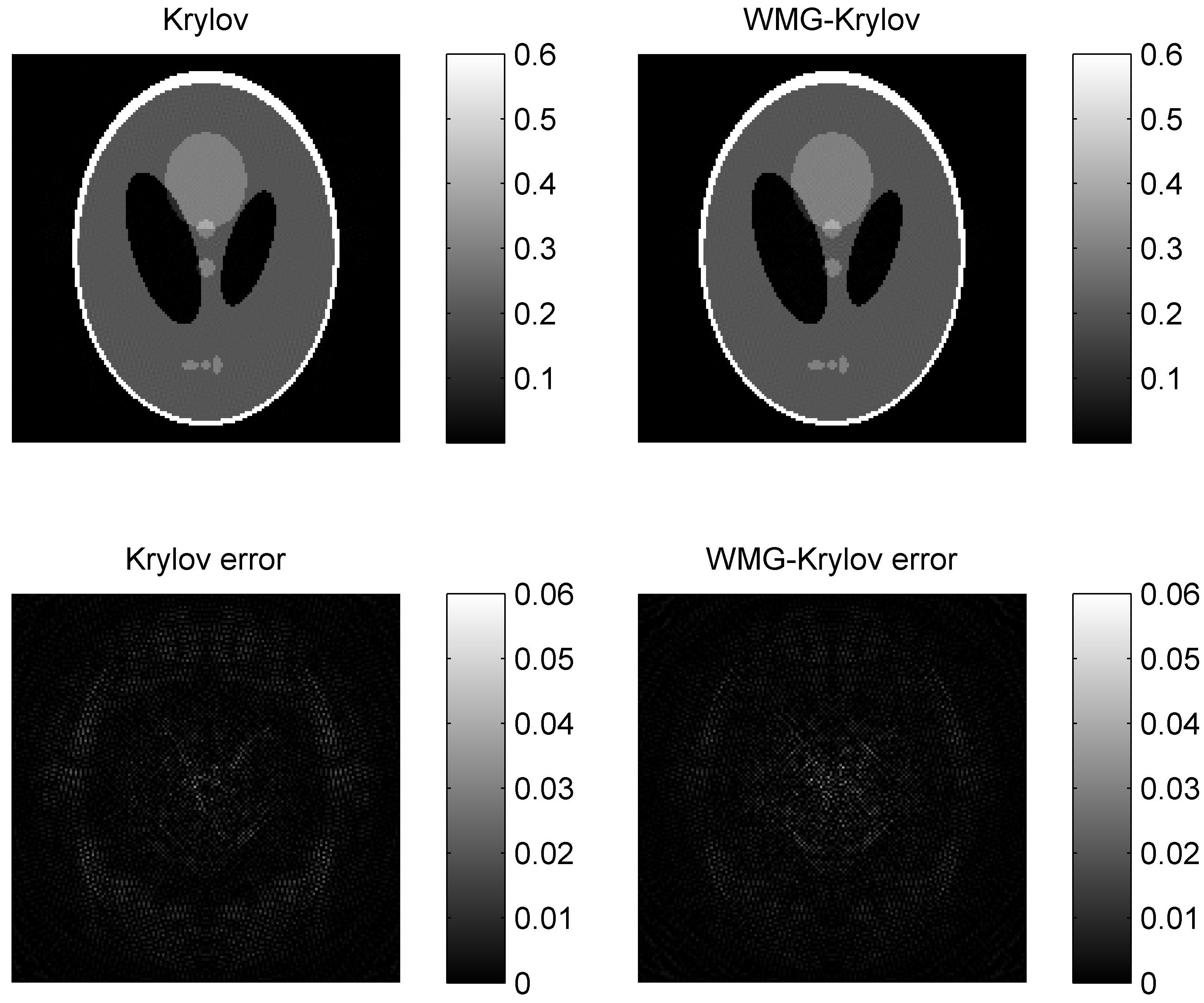}
\caption{Shepp-Logan type model problem with $N = 160 \times 160$ and $M = 400 \times 160$ (no noise, regularized). Displayed are the numerical solutions computed by 300 BiCGStab iterations (left panel) and 50 WMG-preconditioned BiCGStab iterations (right panel), with a small regularization parameter $\lambda = 0.4$ for both methods. Specifications: see Table \ref{tab:resultss}.}
\label{fig:soll}
\end{center}
\end{figure}

\begin{table}[t]
\centering
\begin{tabular}{|c| c c c|}
\hline
  & iterations & CPU time & $L_2$ error \\
\hline
  {\footnotesize BiCGStab} 		 & 300   & 25.5 s. & 0.0180 \\
  {\footnotesize WMG-BiCGStab} & 50    & 17.7 s. & 0.0165 \\
\hline
\end{tabular}
\caption{Shepp-Logan type model problem with $N = 160 \times 160$ and $M = 400 \times 160$ (no noise, regularized). Displayed are the number of iterations required to reach the error criterion (tol = 2\%), elapsed CPU time (in seconds), and relative $L_2$ error $\|x_k - x_{ex}\|_2/\|x_{ex}\|_2$ for different solution methods.}
\label{tab:resultss}
\end{table}

When solving a noisy and/or underdetermined ill-posed system, the discrepancy between the exact non-noisy object of interest $x_{ex}$ and the noisy reconstructed image $x^*$ has to be accounted for. This is commonly done by the inclusion of a regularization term in the linear system (\ref{eq:sys}). In this work, we apply standard Tikhonov regularization \cite{tikhonov1995numerical,golub1999tikhonov}, minimizing the regularized residual
\beq
\min_x  \{ \| Wx - \tilde{b} \|_2 + \lambda \| x \|_2 \}.
\eeq
This is equivalent to solving the regularized system of normal equations
\beq \label{eq:sys_reg}
(W^T W + \lambda) x = W^T \tilde{b},
\eeq
where $\lambda$ is the regularization parameter which is generally chosen to be small with respect to the spectral radius of $W^T W$. We remark that the exact value of $\lambda$ depends on the solution method. 

In addition to the suppression of noise (see Section 4.3), regularization naturally eliminates small-scale artifacts like the one rendered by the WMG preconditioner from the reconstruction. Figure \ref{fig:soll} shows the Krylov and WMG-Krylov solutions to the model problem from Section 4.1, with the inclusion of a relatively small regularization parameter $\lambda = 0.4$ for both methods. The resulting solutions are of comparable quality, see the corresponding Table \ref{tab:resultss}. However, the regularization term removes the artifact in the center of the WMG-Krylov reconstruction, yielding a relative $L_{\infty}$ error norm of $3.99\%$, which is comparable to the Krylov method $L_{\infty}$ error norm ($3.19\%$).

\subsection{Shepp-Logan type model problem with noisy data}

The numerical experiments in the previous section were primarily intended as an academic test case to demonstrate the performance of the WMG preconditioner. In this section we consider a more realistic model problem by incorporating artificially generated noise into the data. In realistic applications, one aims at solving system (\ref{eq:sys})-(\ref{eq:sys_ne}) with a noisy right-hand side, 
\beq \label{eq:sys_noisy}
W x = \tilde{b},
\eeq
where $\tilde{b} = b + \varepsilon$. Here $b$ represents the non-noisy projection data and $\varepsilon$ is a noise term. We consider a noise vector of randomly distributed white noise 
\beq \label{eq:noise_term}
\varepsilon = \alpha \, U(-1,1) \max(|b|),
\eeq
where $\alpha$ is the noise level, commonly given in $\%$ of $\max(|b|)$, and $U(-1,1)$ designates the realization of an $M$-dimensional random variable selected from the $M$-dimensional uniform distribution on the open interval $(-1,1)$. We note that the exact object of interest $x_{ex}$ is generally not a solution to system (\ref{eq:sys_noisy}). Let the solution to (\ref{eq:sys_noisy}) be denoted by $x^*$, then 
\beq 
x^* = x_{ex} + \delta,
\eeq
where $\delta$ is the discrepancy between the exact non-noisy target image $x_{ex}$ and the noisy solution $x^*$ to (\ref{eq:sys_noisy}). This discrepancy $\delta$ is the backprojection of the noise term $\varepsilon$.

\begin{figure}[t]
\begin{center}
\includegraphics[width=0.55\textwidth]{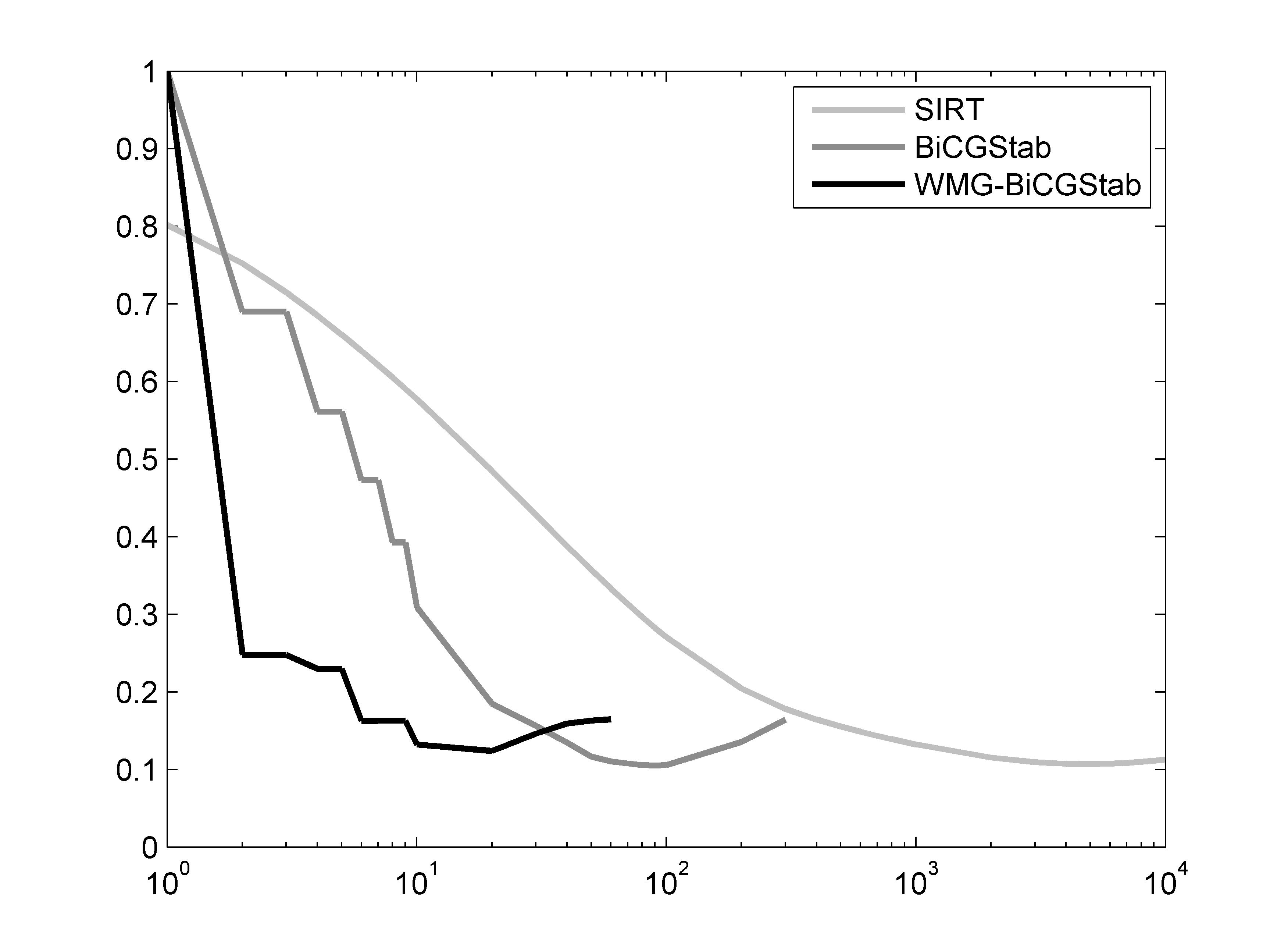}
\caption{Shepp-Logan type model problem with $N = 160 \times 160$ and $M = 400 \times 160$ (noisy, non-regularized). Displayed are the scaled error $L_2$ norms $\|e_k\|_2/\|e_0\|_2$ for the SIRT (light gray), BiCGStab (dark gray) and WMG-preconditioned BiCGStab (black) method in function of the number of iterations. Error minima at $k_{opt} = 3000$ (SIRT), $k_{opt} = 100$ (BiCGStab) and $k_{opt} = 14$ (WMG-BiCGStab) iterations. Horizontal axis in log scale.}
\label{fig:semiconv}
\end{center}
\end{figure}

\subsubsection{Semi-convergence.}

In the inverse problems literature \cite{hansen2010discrete,engl1996regularization,siltanen2012linear}, a typical convergence phenomenon referred to as \emph{semi-convergence} is described, which occurs when trying to recover the exact solution $x_{ex}$ from a noisy system like (\ref{eq:sys_noisy}) using iterative methods. Comparing the iterative solution to (\ref{eq:sys_noisy}) with the exact image $x_{ex}$ in each step, we observe that convergence of all methods tends to stagnate gradually, with the error $\|e_k\|_2$ reaching a minimum after a certain problem-dependent (up-front unknown) number of iterations. Beyond this point the error increases, and additional iterations push the iterative solution increasingly further away from the exact solution. Note that this behaviour is not apparent from the residual history, as the residual tends to keep decreasing with every iteration. For each iterative method, there exists an optimal number of iterations $k_{opt}$ which minimizes the error,
\beq
	k_{opt} = \underset{k}{\operatorname{arg\,min}} \, \|e_k\|_2.
\eeq
To illustrate this concept, we recall the 160-by-160 pixels Shepp-Logan type model problem with $400 \times 160$ element data vector, as introduced in Section 4.1, and we artificially add $\alpha = 1\%$ of white noise to the data vector $b$. Figure \ref{fig:semiconv} shows the scaled $L_2$ norm of the error $\|e_k\|_2/\|e_0\|_2$ for the SIRT method, the unpreconditioned Krylov method and the WMG-Krylov method for the noisy model problem. As predicted by the analysis in \cite{hansen2010discrete}, an error minimum is reached for all methods, with SIRT reaching a minimum after 3000 iterations, while the Krylov and WMG-Krylov method require 100 and 14 iterations respectively to reach the error minimum.

\subsubsection{Regularizing the noisy problem.}

\begin{figure}[t]
\begin{center}
\includegraphics[width=0.55\textwidth]{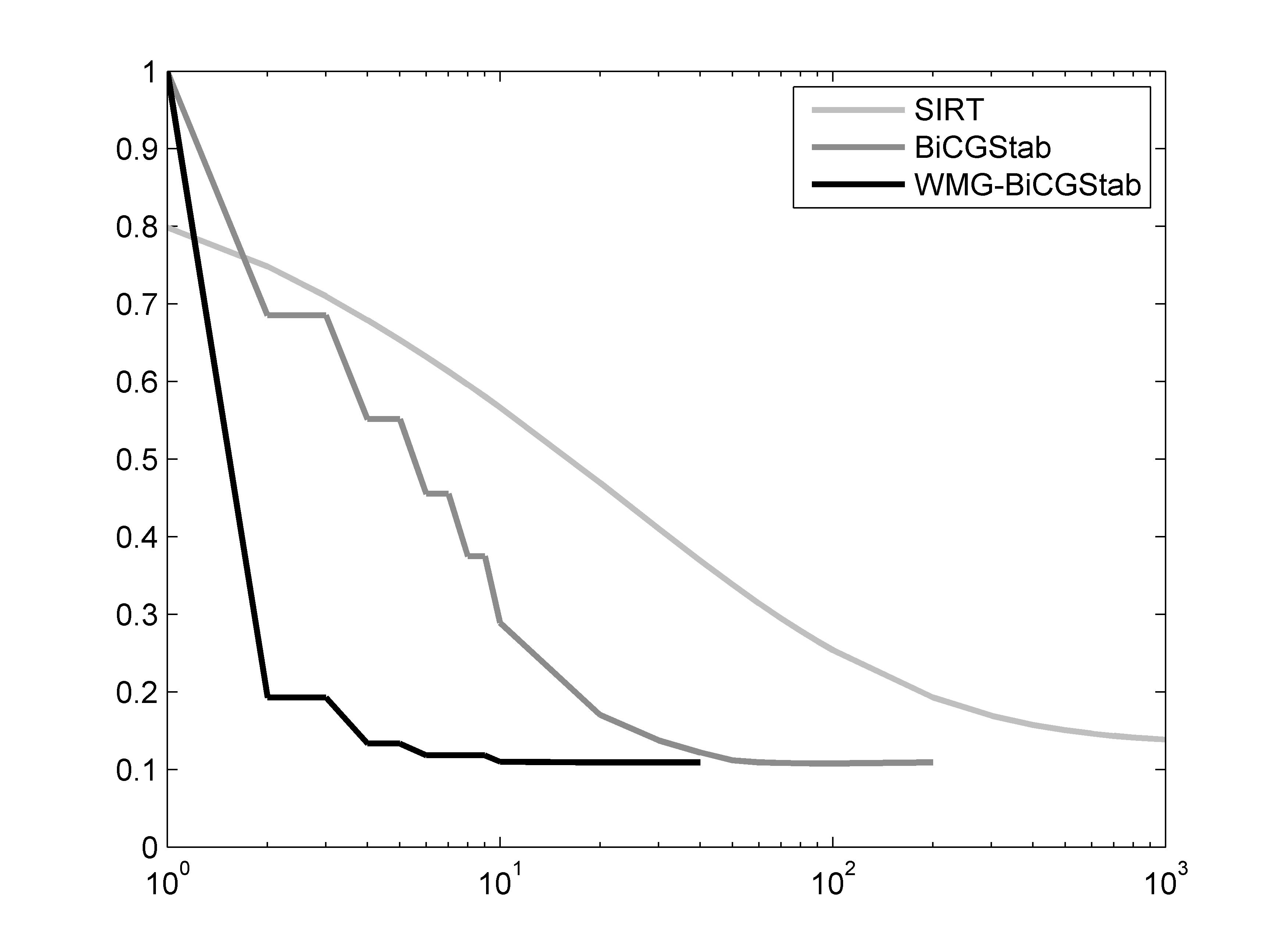}
\caption{Shepp-Logan type model problem with $N = 160 \times 160$ and $M = 400 \times 160$ (noisy, regularized). Displayed are the scaled error $L_2$ norms $\|e_k\|_2/\|e_0\|_2$ for the SIRT (light gray), BiCGStab (dark gray), and WMG-preconditioned BiCGStab (black) method in function of the number of iterations. Error minima at $k_{opt} = 3000$ (SIRT), $k_{opt} = 100$ (BiCGStab), and $k_{opt} = 14$ (WMG-BiCGStab) iterations. Regularization parameter: SIRT: $\lambda = 0.001$, Krylov: $\lambda = 10$, WMG-Krylov: $\lambda = 10$. Horizontal axis in log scale.}
\label{fig:semiconv_relax}
\end{center}
\end{figure}

As proposed in standard works on regularization of inverse problems \cite{hansen2010discrete, engl1996regularization}, Tikhonov regularization is often imposed to ensure proper convergence of the iterative solvers towards the non-noisy solution.
Note that Figure \ref{fig:semiconv} proves the existence of an error minimum in the non-regularized case, but does not predict a realistic estimation of the regularized error. Figure \ref{fig:semiconv_relax} is the analogon of Figure \ref{fig:semiconv} with the inclusion of a Tikhonov regularization term.
When the regularization parameter $\lambda$ is well-chosen, the regularization term ensures a denoising (or smoothing) of the reconstruction in every iteration. Consequently, the error $L_2$ norm $\|e_k\|_2$ decreases in every iteration, even after reaching the optimal number of iterations $k_{opt}$, see Figure \ref{fig:semiconv_relax}. As shown by Hansen in \cite{hansen2010discrete}, the error decreases only marginally from this point onwards and the iteration is typically stopped after $k_{opt}$ iterations.

\subsubsection{Performance results.}

\begin{figure}[t]
\begin{center}
\includegraphics[width=1\textwidth]{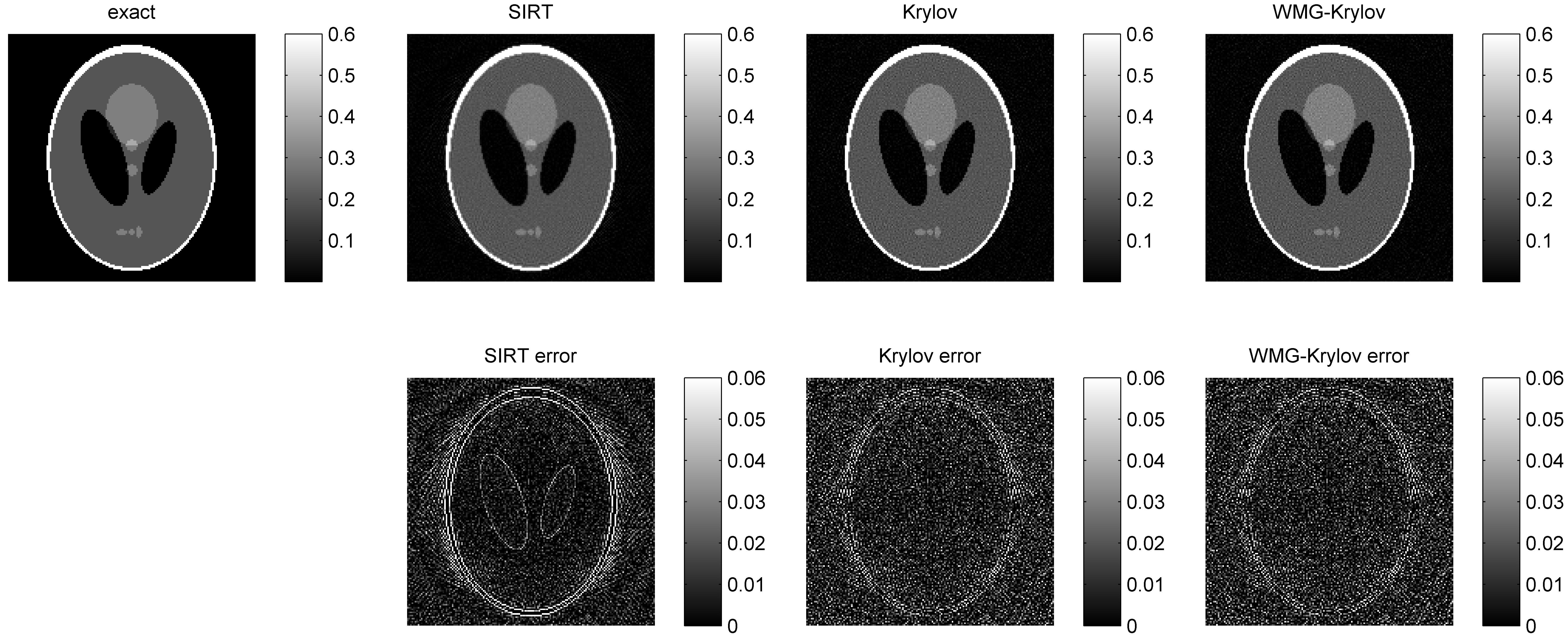}
\caption{Shepp-Logan type model problem with $N = 160 \times 160$ and $M = 400 \times 160$ (noisy, regularized). Displayed are the numerical solutions to the regularized system (\ref{eq:sys_reg}) after $k_{opt} = 1000$ (SIRT), $k_{opt} = 100$ (BiCGStab) and $k_{opt} = 14$ (WMG-BiCGStab) iterations. Regularization parameters: SIRT: $\lambda = 0.001$, Krylov: $\lambda = 10$, WMG-Krylov: $\lambda = 10$. Specifications: see Table \ref{tab:results2}.}
\label{fig:sol_noise}
\end{center}
\end{figure}

We again consider the 160-by-160 pixels Shepp-Logan type model problem with the addition of $\alpha = 1\%$ white noise. Figure \ref{fig:sol_noise} compares the resulting reconstructions after $k_{opt} = 1000$ SIRT iterations, $k_{opt} = 100$ BiCGStab iterations and $k_{opt} = 14$ WMG-BiCGStab iterations respectively. The regularization parameters are chosen as $\lambda = 0.001$ for the SIRT iteration and $\lambda = 10$ for the Krylov methods. Note that the regularization parameter for the SIRT method is smaller than the Krylov parameter due to rescaling of the SIRT system (\ref{eq:sys_sirt}). 
An overview of the corresponding reconstruction characteristics can be found in Table \ref{tab:results2}. Note that the accuracy of the Krylov reconstructions is comparable; however, the Krylov solutions feature less artifacts and a sharper characterization of the edges, resulting in a smaller $L_2$ and $L_{\infty}$ norm compared to the SIRT method. Additionally, the WMG-BiCGStab method uses only 14 iterations to reach an accuracy comparable to that of the BiCGStab solution after 100 iterations. This results in a computational time of less than 5.5 seconds by the WMG-BiCGStab method, which is significantly lower than the time elapsed by the BiCGStab iterations (9.35 s.)~in order to obtain a comparable accuracy. A CPU time speed-up of approximately 42\% is achieved by the WMG preconditioner. The improved convergence speed of the WMG-Krylov method is a major advantage over classical Algebraic Reconstruction Techniques, and the WMG preconditioner features a significantly reduced number of iterations compared to standard unpreconditioned Krylov solvers.

\begin{table}[t]
\centering
\begin{tabular}{|c| c c c c|}
\hline
   & iterations & CPU time & error $L_2$ & error $L_{\infty}$ \\
\hline
  {\footnotesize SIRT} 				 & 1000 & 81.1 s. & 0.1385 & 0.2697 \\
  {\footnotesize BiCGStab} 		 & 100  & 9.35 s. & 0.1074 & 0.1459 \\
  {\footnotesize WMG-BiCGStab} & 14   & 5.44 s. & 0.1083 & 0.1386 \\
\hline
\end{tabular}
\caption{Shepp-Logan type model problem with $N = 160 \times 160$ and $M = 400 \times 160$ (noisy, regularized). Displayed are the number of iterations $k_{opt}$ (see Figure \ref{fig:semiconv}), elapsed CPU time (in seconds), relative $L_2$ error $\|e_k\|_2/\|e_0\|_2$, and relative $L_{\infty}$ error $\|e_k\|_{\infty}/\|e_0\|_{\infty}$ for different solution methods.}
\label{tab:results2}
\end{table}

\section{Conclusions and discussion}
In this paper we proposed a novel algebraic reconstruction method for the linear inversion problems that arise from computerized tomographic reconstruction. Driven by the observed slow convergence of classical Algebraic Reconstruction Techniques, Krylov methods are suggested as a more efficient alternative for algebraic tomographic reconstruction. However, these methods are known to be generally uncompetitive without a suitable preconditioner.

Inspired by the work done by McCormick et al.~\cite{mccormick1993multigrid,henson1996multilevel} and R\"ude et al.~\cite{kostler2006towards}, a multi-level type preconditioning approach is suggested. An eigenvalue analysis of the SIRT method and the classical multigrid scheme shows that standard multigrid is unsuitable as a Krylov preconditioner for algebraic tomographic reconstruction. Consequently, a novel wavelet-based multigrid (WMG) preconditioner is introduced that projects the large fine-level matrix operator onto a collection of smaller coarse level subproblems. The advantage of this approach is that the coarse grid subproblems are computationally cheaper to solve, resulting in a fast and efficient overall preconditioning scheme. It is shown that the WMG-preconditioned Krylov method has improved spectral properties, and thus yields a performant iterative solver for tomographic reconstruction.

Additionally, when the domain and the matrix operator $W$ become larger and no longer fit in fast cache memory, performance of the SpMV, and thus stationary iterations like SIRT, typically decreases. This is because the SpMV is memory bandwidth limited, cf.~\cite{vanAarle2012memory}. The WMG preconditioner suffers much less from this effect, because it applies (small) dense matrix operations (triangular solves) on the coarsest level, which achieve high efficiency because they have a high flop-to-memory-access ratio and are implemented in highly optimized BLAS3 routines. Due to the strong reduction in the number of Krylov iterations and the scalable nature of the WMG scheme, the WMG-preconditioned Krylov method is expected to be increasingly performant when solving large-scale and 3D reconstruction problems. 

The WMG-preconditioned Krylov method is primarily analyzed on a 2D non-noisy benchmark problem, where a notable convergence speedup is observed compared to standard unpreconditioned Krylov methods. Furthermore, the approach is successfully validated on a more realistic noisy benchmark problem, where Tikhonov regularization is introduced to account for the ill-posedness of the problem. 
Comparing the new WMG-Krylov method to existing Algebraic Reconstruction Techniques like the classical SIRT scheme and unpreconditioned Krylov methods, we find a comparable quality of the reconstruction at a significantly reduced number of iterations. Numerical experiments on semi-realistic benchmark problems confirm that the WMG-Krylov method outperforms the classical reconstruction methods in terms of computational time, while retaining a comparable precision on the numerical solution.

We note that the WMG method can in principle be applied directly to the (possibly underdetermined) system (\ref{eq:sys}), cf.~the work done by R\"ude et al.~in \cite{kostler2006towards}. The embedding in a governing Krylov solver generally improves robustness, but requires the use of the normal equations formulation (\ref{eq:sys_ne}). If the WMG method is applied directly to the underdetermined system, restriction can theoretically be applied on both sides of the system, i.e.~simultaneously restricting the number of pixels and the number of rays or projection angles. This strategy could be explored in future research. 

To conclude, the work done in this paper is an initial effort to designing a performance based preconditioning technique for algebraic tomographic reconstruction problems. Whereas preconditioning approaches suggested in other works focus primarily on regularization properties (noise reduction, deblurring, etc.), the most important feature of the WMG preconditioner is the significant reduction in the number of Krylov iterations, which allows to effectively speed-up Krylov convergence, and hence reduces the reconstruction time. 

\section{Acknowledgments}
\noindent This research was partly funded by the \textit{Fonds voor Wetenschappelijk Onderzoek (FWO)} project G.0.120.08 and \textit{Krediet aan navorser} project number 1.5.145.10. Additionally, this work was partly funded by Intel and the \textit{Institute for the Promotion of Innovation through Science and Technology in Flanders (IWT)}.

\section*{References}
\bibliography{refs}

\end{document}